\documentclass[11pt]{amsart}
\usepackage{amsmath}
\usepackage{amssymb}
\usepackage{amsfonts}
\usepackage{amscd}
\usepackage[mathscr]{eucal}
\usepackage{graphicx}
\usepackage{float}

\def\ds{\displaystyle}

\newtheorem{Theorem}{Theorem}[section]
\newtheorem{Corollary}[Theorem]{Corollary}
\newtheorem{Lemma}[Theorem]{Lemma}
\newtheorem{Proposition}[Theorem]{Proposition}

\theoremstyle{definition}
\newtheorem{Definition}[Theorem]{Definition}
\newtheorem{Example}[Theorem]{Example}
\newtheorem{Question}[Theorem]{Question}
\newtheorem{Problem}[Theorem]{Problem}

\theoremstyle{remark}
\newtheorem{Remark}[Theorem]{Remark}

\title[On rigidity of Coxeter systems up to finite twists]{
On rigidity of Coxeter systems up to finite twists and 
separations of Coxeter generating sets}
\author{Tetsuya Hosaka}
\address{Department of Mathematics, Faculty of Science, Shizuoka University, 
836 Ohya, Suruga-ku, Shizuoka, 422-8529, Japan}
\date{June 20, 2023}
\email{hosaka.tetsuya@shizuoka.ac.jp}
\thanks{This work was partly supported by JSPS KAKENHI Grant Number JP18K03273.}

\begin{document}

\sloppy

\begin{abstract}
In this paper, 
we study the twist-conjecture for Coxeter systems and rigidity of Coxeter systems up to finite twists.
For Coxeter systems $(W,R)$ and $(W,S)$, 
under the untangle-condition for conjugate subsets, 
we investigate separations and type(I) and type(II) subsets of $R$ and $S$ 
and give an equivalent condition of $R$ and $S$ 
that are conjugate up to finite twists.
We provide one direction of approach 
to solving the twist-conjecture and the isomorphism problem for Coxeter groups of finite ranks.
\end{abstract}

\maketitle

\section{Introduction and preliminaries}\label{sec1}
The purpose of this paper is to 
investigate the twist-conjecture for Coxeter systems and rigidity of Coxeter systems up to finite twists.
Definitions and details of Coxeter groups and Coxeter systems 
are found in \cite{Bo} and \cite{Hu}.
In this paper, we suppose that all Coxeter groups are of finite ranks.
The isomorphism problem for Coxeter groups is open.

\begin{Problem}[The isomorphism problem for Coxeter groups]
For given Coxeter systems $(W,S)$ and $(W',S')$, 
find an algorithm to determine whether 
the Coxeter groups $W$ and $W'$ are isomorphic or not.
\end{Problem}

For a Coxeter system $(W,S)$ and a subset $T$ of $S$, 
$W_T$ is defined as the subgroup generated by $T$ in $W$.
It is well known that $(W_T,T)$ also becomes a Coxeter system.
A subset $T$ of $S$ is said to be \textit{spherical}, 
if the subgroup $W_T$ is finite.
If $T$ is a spherical subset of $S$, then 
it is known that there exists a unique longest length element 
in $W_T$ that is denoted by $w_T$.
For a subset $T$ of $S$, 
we denote $T^\perp:=\{s\in S : st=ts \ \text{for all}\ t\in T\}$.
For a Coxeter group $W$, 
a subset $R$ of $W$ is called a \textit{Coxeter generating set} for $W$, 
if $(W,R)$ is a Coxeter system.

Recent research on the isomorphism problem for Coxeter groups 
and rigidity of Coxeter groups and Coxeter systems can be found 
in \cite{BMMN}, \cite{CM}, \cite{CP}, \cite{HP}, 
\cite{MM}, \cite{Mu}, \cite{MuWe}, \cite{N}, \cite{RT} and \cite{W}.

A Coxeter system $(W,S)$ is said to be \textit{reflection-rigid}, 
if it is determined by the Coxeter group $W$ and the set of reflections 
\[ {\mathcal R}_S:=\{wsw^{-1}\in W : w\in W,\ s\in S\} =S^W \]
up to isomorphisms; 
that is, 
for any Coxeter generating set $S'$ for $W$ such that 
the reflection sets ${\mathcal R}_S$ and ${\mathcal R}_{S'}$ are equal, 
two Coxeter systems $(W,S)$ and $(W,S')$ are isomorphic.

Also a Coxeter system $(W,S)$ is said to be \textit{strongly-reflection-rigid},
if $S$ is determined by the Coxeter group $W$ and the reflection set 
${\mathcal R}_S$ up to conjugate; 
that is, 
for any Coxeter generating set $S'$ for $W$ such that 
the reflection sets ${\mathcal R}_S$ and ${\mathcal R}_{S'}$ are equal, 
two Coxeter generating sets $S$ and $S'$ are conjugate in $W$.

Let $W$ be a Coxeter group.
Two Coxeter generating sets $S$ and $S'$ for $W$ are said to be 
\textit{angle-compatible} \cite{CP} (or \textit{sharp-angled} \cite{MM}), if 
for any $\{a,b\}\subset S$ as $o(ab)<\infty$ 
there exists $\{a',b'\}\subset S'$ as $o(a'b')<\infty$ 
such that $\{a,b\}$ and $\{a',b'\}$ are conjugate in $W$.
Here for $w\in W$, $o(w)$ is the order of $w$ in the group $W$.
We say that a Coxeter system $(W,S)$ is \textit{angle-rigid}, if 
for any Coxeter generating set $S'$ for $W$ such that 
$S$ and $S'$ are angle-compatible, 
two Coxeter systems $(W,S)$ and $(W,S')$ are isomorphic.
Also a Coxeter system $(W,S)$ is said to be \textit{strongly-angle-rigid}, 
if for any Coxeter generating set $S'$ for $W$ such that 
$S$ and $S'$ are angle-compatible, 
two Coxeter generating sets $S$ and $S'$ are conjugate in $W$.

It is known that 
if two Coxeter generating sets $S$ and $S'$ for $W$ are angle-compatible 
then the reflection sets ${\mathcal R}_S$ and ${\mathcal R}_{S'}$ 
are equal (see \cite{CP}).
Hence 
if a Coxeter system $(W,S)$ is (strongly-)reflection-rigid, 
then it is (strongly-)angle-rigid.

Recently, Marquis and M\"{u}hlherr have proved the following.

\begin{Theorem}[{\cite[Corollary~1.1]{MM}}]\label{thmMM}
The isomorphism problem for Coxeter groups is solved 
as soon as the following problem is solved.
\end{Theorem}

\begin{Problem}[{\cite{MM}}]\label{problem1.2}
Let $(W,R)$ be a Coxeter system.
Find all Coxeter generating sets $S \subset R^W$ 
such that $S$ is sharp-angled with respect to $R$.
\end{Problem}

Thus if for a given Coxeter system $(W,R)$ 
we can output all Coxeter generating sets $S$ for $W$ such that 
$R$ and $S$ are angle-compatible, 
then the isomorphism problem for Coxeter groups is solved.

Caprace and Przytycki have given the following.

\begin{Theorem}[{\cite[Theorem~1.1]{CP}}]\label{thmCP}
Let $S$ and $R$ be angle-compatible Coxeter generating sets 
for a group $W$.
If $S$ is twist-rigid, then $S$ and $R$ are conjugate.
\end{Theorem}

Hence, every twist-rigid Coxeter system is strongly-angle-rigid.

Here $(W,S)$ (and $S$) is said to be \textit{twist-rigid}, if 
$(W,S)$ has no elementary twist \cite{CP}.
We give a definition and detail of a ``twist'' later.

\begin{Definition}\label{def:essential}
Let $(W,S)$ be a Coxeter system and 
let $U$ be a subset of $S$.
We consider the irreducible decomposition 
\[ W_U= W_{U_1} \times \cdots \times W_{U_n}, \]
where $U=U_1 \cup \cdots \cup U_n$ is a disjoint union and 
each $(W_{U_i},U_i)$ is an irreducible Coxeter system.
Here each $W_{U_i}$ is either finite or infinite.
We define 
\begin{align*} 
&U_\sigma := \bigcup \{U_i : \text{$W_{U_i}$ is finite} \} \ \text{and} \\ 
&U_\nu := \bigcup \{U_i : \text{$W_{U_i}$ is infinite} \}.
\end{align*}
Then $U=U_\sigma \cup U_\nu$ is a disjoint union, 
$W_{U_\sigma}$ is finite, $W_{U_\nu}$ is infinite (if $U_\nu$ is non-empty), 
and $W_U=W_{U_\sigma} \times W_{U_\nu}$.
\end{Definition}

Let $(W,S)$ be a Coxeter system.
We say that a subset $U$ of $S$ is 
a \textit{spherical-product subset} of $S$, 
if $U$ is non-empty and 
$U \subset \sigma \cup \sigma^\perp$ 
for some non-empty spherical subset $\sigma$ of $S$.
Here if $U$ is a spherical-product subset of $S$, then 
either 
\begin{enumerate}
\item $U_\sigma =U$ and $U_\nu =\emptyset$ (where $U$ is spherical), 
\item $U_\sigma \neq \emptyset$ and $U_\nu \neq \emptyset$ 
(where $U=U_\sigma \cup U_\nu \subset U_\sigma \cup (U_\sigma)^\perp$), or
\item $U_\sigma =\emptyset$ and $U_\nu =U$ 
(where for a spherical subset $\sigma$ of $S$ 
such that $U \subset \sigma \cup \sigma^\perp$, 
if we put $\overline{U}:=\sigma \cup U$ then 
$\overline{U}_\sigma=\sigma$ and $\overline{U}_\nu=U$).
\end{enumerate}

A subset $T$ of $S$ is said to be \textit{connected}, 
if for any $a,b \in T$ as $a\neq b$, 
there exists a sequence $a=t_1,t_2,\ldots,t_n=b$ in $T$ 
such that $o(t_i t_{i+1})$ is finite for any $i=1,\ldots,n-1$ 
(that is, the nerve of $(W_T,T)$ is connected).
Also we often say that a spherical-product subset $U$ 
\textit{separates} $S$, 
if $S$ is connected and $S-U$ is not connected 
(that is, the nerve $L_U$ of $(W_U,U)$ separates 
the connected nerve $L$ of $(W,S)$ to 
some at least two components).

If a spherical-product subset $U$ of $S$ separates $S$, 
then for some non-empty subsets $X$ and $Y$ of $S$, 
\begin{enumerate}
\item[(1)] $S-U=X \cup Y$ that is a disjoint union and 
\item[(2)] $o(xy)=\infty$ for any $x\in X$ and $y\in Y$.
\end{enumerate}
There is a spherical subset $\sigma$ of $S$ such that 
$U \subset \sigma \cup \sigma^\perp$.
Here if $U_\sigma$ is non-empty then we take $\sigma:=U_\sigma$.
Then $w_\sigma U w_\sigma=U$ and 
we obtain a Coxeter generating set for $W$ as 
\[ S':= X \cup U \cup (w_\sigma Y w_\sigma) \]
that is an \textit{elementary twist} (\cite{BMMN}, \cite{BH}, \cite{CP}, \cite{Mu0}, \cite{Mu}). 
Also in the case that $S$ is not connected, 
we can consider that $U=\emptyset$ separates $S$.
If $S=X \cup Y$ that is a disjoint union and 
$o(xy)=\infty$ for any $x\in X$ and $y\in Y$, then 
for any $w\in W_X$ we obtain a Coxeter generating set 
$S':= X \cup (w Y w^{-1})$ for $W$ 
that is also an elementary twist (\cite{BMMN}, \cite{BH}, \cite{CP}, \cite{CD}, \cite{Mu}). 

Hence, $(W,S)$ is twist-rigid if and only if $S$ is connected and 
any spherical-product subset $U$ does not separate $S$.

More generally, let $U$ be 
a spherical-product subset of $S$ that separates $S$ and 
let $X$ and $Y$ be non-empty subsets of $S$ such that 
\begin{enumerate}
\item[(1)] $S-U=X \cup Y$ that is a disjoint union, 
\item[(2)] $o(xy)=\infty$ for any $x\in X$ and $y\in Y$ and
\item[(3)] $wUw^{-1} =U$ for some $w\in W_{X\cup U}$.
\end{enumerate}
Then we obtain a Coxeter generating set for $W$ as
\[ S':= X \cup U \cup (w Y w^{-1}) \]
that is called a \textit{$($general-$)$twist}.
Here there is a possibility that 
every twist can be denoted by some elementary-twists.
The author does not know whether this always holds.

In this paper, we often say that 
a Coxeter generating set $S'$ is obtained from $S$ 
by some finite twists, 
if there exists a sequence 
$S=S_1,S_2,\ldots,S_n=S'$ of Coxeter generating sets for $W$ 
such that each $S_{i+1}$ is obtained from $S_i$ 
by some twist.
We also say that 
Coxeter generating sets $R$ and $S$ for $W$ 
are \textit{conjugate up to finite twists}, if 
there exists a Coxeter generating set $R'$ for $W$ 
such that $R'$ is obtained from $R$ by some finite twists and 
$R'$ is conjugate to $S$.
Here Coxeter generating sets $R$ and $S$ for $W$ 
are conjugate up to finite twists if and only if 
there exist Coxeter generating sets $R'$ and $S'$ for $W$ 
such that $R'$ and $S'$ are obtained from $R$ and $S$ 
by some finite twists respectively 
and $R'$ is conjugate to $S'$.

Now we define new concepts 
``separations'' and ``untangle-conjugate'' on Coxeter systems.
These definitions are technical.
The purpose of this paper is to give a further reduction of the isomorphism problem for Coxeter groups.
These definitions are designed and constructed as the main results in this paper hold and 
as two Coxeter generating sets with some conditions using them are conjugate up to finite twists.

First, the following lemma is known.

\begin{Lemma}[cf.\ \cite{BMMN}]\label{lemma0-0}
Let $(W,S)$ and $(W,S')$ be Coxeter systems.
For each maximal spherical subset $T$ of $S$, 
there exists a unique maximal spherical subset $T'$ of $S'$ 
such that $W_T$ and $W_{T'}$ are conjugate in $W$.
\end{Lemma}

Here we say that 
two Coxeter systems $(W,S)$ and $(W,S')$ 
are \textit{maximal-spherical-subset-compatible}, 
if for each maximal spherical subset $T$ of $S$, 
there exists a maximal spherical subset $T'$ of $S'$ 
such that $T$ and $T'$ are conjugate in $W$ 
(where $T'$ is unique by Lemma~\ref{lemma0-0}).
If two Coxeter systems are 
maximal-spherical-subset-compatible 
then they are angle-compatible. 
Indeed for each $\{a,b\} \subset S$ as $o(ab)<\infty$, 
there exists a maximal spherical subset $T$ of $S$ 
containing $\{a,b\}$.
If $T'$ is a maximal spherical subset of $S'$ 
such that $T$ and $T'$ are conjugate in $W$ 
then some subset $\{a',b'\} \subset T'$ is conjugate 
to $\{a,b\}$.

For a Coxeter system $(W,S)$ as $S$ is connected 
and for $A\subset S$, 
we say that $A$ is a \textit{twist-rigid subset} of $S$, 
if $A$ is connected and 
if there does not exist a spherical-product subset $U$ of $S$ 
such that $U$ separates $S$ and $U\cap A$ separates $A$.
Here we note that $(W_A,A)$ need not be a twist-rigid Coxeter system 
in general (see examples in Section~\ref{sec2}).

Let $(W,S)$ be a Coxeter system.
We suppose that $S$ is connected.
Let ${\mathcal A}$ be a set of subsets of $S$ such that 
\begin{enumerate}
\item[(1)] $S=\bigcup_{A \in \mathcal A} A$,
\item[(2)] each $A\in {\mathcal A}$ is connected and 
a union of some maximal twist-rigid subsets of $S$, and
\item[(3)] for each maximal twist-rigid subset $A_0$ of $S$, 
there exists a unique element $A\in {\mathcal A}$ such that $A_0 \subset A$.
\end{enumerate}
A subset $U$ of $S$ is called a \textit{separator} of ${\mathcal A}$, 
if the following conditions (i)--(v) hold:
\begin{enumerate}
\item[(i)] $U$ is a spherical-product subset of $S$.
\item[(ii)] $U$ separates $S$.
\item[(iii)] For any $A\in {\mathcal A}$, 
there exists a unique $j \in \{1,\ldots,t\}$ such that $A\subset \overline{X}_j$.
Here $S-U= X_1 \cup \cdots \cup X_t$ and 
$X_1, \ldots, X_t$ are the connected components of $S-U$ and 
we define $\overline{X}_j:=\bigcup \{ A\in {\mathcal A} : A \subset X_j\cup U \}$ 
for each $j=1,\ldots,t$.
\item[(iv)] There exist $A_1,A_2\in {\mathcal A}$ such that $A_1 \cap A_2 =U$, 
$A_1 \subset \overline{X}_{j_1}$ and $A_2 \subset \overline{X}_{j_2}$ 
for some $j_1, j_2\in \{1,\ldots,t\}$ as $j_1\neq j_2$.
\item[(v)] Let $j\in \{1,\ldots,t\}$.
For any distinct elements $A_1,\ldots,A_n \in {\mathcal A}$ such that 
$U \cap \overline{X}_j \subset A_1 \cup \cdots \cup A_n\subset \overline{X}_j$, 
if $(A_1\cup \cdots \cup A_i) \cap A_{i+1}$ is 
maximal in 
\[ \{(A_1\cup \cdots \cup A_i) \cap A : A \in {\mathcal A}-\{A_1,\ldots,A_i\} \} \]
for any $i=1,\ldots,n-1$, 
then $U \cap \overline{X}_j \subset A_i$ for some $i\in \{1,\ldots,n\}$.
\end{enumerate}

A set ${\mathcal A}$ of subsets of $S$ is called 
a \textit{separation} of $S$, 
if the following conditions (1)--(4) hold:
\begin{enumerate}
\item[(1)] $S=\bigcup_{A \in \mathcal A} A$.
\item[(2)] Each $A\in {\mathcal A}$ is connected and 
a union of some maximal twist-rigid subsets of $S$.
\item[(3)] For each maximal twist-rigid subset $A_0$ of $S$, 
there exists a unique element $A\in {\mathcal A}$ such that $A_0 \subset A$.
\item[(4)] 
For distinct elements $A_1,A_2,\ldots,A_n \in {\mathcal A}$, 
if $(A_1\cup \cdots \cup A_i) \cap A_{i+1}$ is 
maximal in 
\[ \{(A_1\cup \cdots \cup A_i) \cap A : A \in {\mathcal A}-\{A_1,\ldots,A_i\} \} \]
for any $i=1,\ldots,n-1$, 
then for each $i=1,\ldots,n-1$, 
\begin{enumerate}
\item[(a)] $U_i:=(A_1\cup \cdots \cup A_i) \cap A_{i+1}$ is a separator of ${\mathcal A}$, and
\item[(b)] $A_1\cup \cdots \cup A_i \subset X_{j_1} \cup U_i$ and $A_{i+1} \subset X_{j_2} \cup U_i$ 
for some $j_1, j_2\in \{1,\ldots,t \}$ as $j_1 \neq j_2$, 
where $S-U_i= X_1 \cup \cdots \cup X_t$ and 
$X_1, \ldots, X_t$ are the connected components of $S-U_i$.
\end{enumerate}
\end{enumerate}

In the case that $S$ is not connected, 
a set ${\mathcal A}$ of subsets of $S$ is called 
a \textit{separation} of $S$, 
if for the connected components $S_1, \ldots,S_k$ of $S$ 
(that is, each $S_i$ is a maximal connected subset of $S$ 
and $S=S_1\cup \cdots \cup S_k$ is a disjoint union), 
${\mathcal A}_i:=\{ A \in {\mathcal A} : A \subset S_i \}$ 
is a separation of $S_i$ for any $i=1,\ldots,k$.

For two sets ${\mathcal A}_1$ and ${\mathcal A}_2$ of subsets of $S$, 
we denote ${\mathcal A}_1 \preceq {\mathcal A}_2$, 
if for any $A_1 \in {\mathcal A}_1$, 
$A_1\subset A_2$ for some $A_2 \in {\mathcal A}_2$.

A separation ${\mathcal A}$ of $S$ is said to be \textit{minimal}, 
if 
\begin{enumerate}
\item[~] there does not exist a separation ${\mathcal A}'$ of $S$ 
such that ${\mathcal A}'\neq {\mathcal A}$ and 
any $A'\in {\mathcal A}'$ is contained in some $A \in {\mathcal A}$.
\end{enumerate}
A minimal separation is ``minimal'' with respect to the partial order ``$\preceq$''.

Let $(W,S)$ be a Coxeter system and 
let ${\mathcal A}$ be a set of subsets of $S$.
We consider a twist 
\[ S':= X \cup U \cup (w Y w^{-1}) \]
where 
$U$ is a spherical-product subset of $S$ and $w\in W_{X \cup U}$ such that 
$wUw^{-1} =U$, 
$U$ separates $S$ as $S-U=X \cup Y$ 
that is a disjoint union and
$o(xy)=\infty$ for any $x\in X$ and $y\in Y$.
A twist $S'$ is said to be \textit{preserving} ${\mathcal A}$, 
if for any $A \in {\mathcal A}$, 
$A \subset X\cup U$ or $A \subset Y\cup U$ 
(that is, there does not exist $A \in {\mathcal A}$ such that 
$A\cap X \neq \emptyset$ and $A\cap Y \neq \emptyset$).
If the twist $S'$ is preserving ${\mathcal A}$ then 
\[ {\mathcal A}':= 
\{ A : A\in {\mathcal A}, \ A \subset X\cup U \} 
\cup \{w A w^{-1} : A\in {\mathcal A},\ A \not\subset X\cup U \} \]
is a set of subsets of $S'$ and 
${\mathcal A}'$ is called 
the set \textit{induced by ${\mathcal A}$ and the twist}.

Here if ${\mathcal A}$ is a separation of $S$, then 
the induced set ${\mathcal A}'$ is a separation of $S'$ 
and it is called 
the \textit{separation induced by ${\mathcal A}$ and the twist}.

Let ${\mathcal A}$ be a set of subsets of $S$.
We often say that 
$S''$ is obtained from $S$ by some finite twists 
\textit{preserving} ${\mathcal A}$, 
if there exists a sequence $S=S_1,S_2,\ldots,S_n=S''$ 
of Coxeter generating sets for $W$ 
such that $S_{i+1}$ is a twist of $S_i$ preserving 
${\mathcal A}_i$ for each $i=1,\ldots, n-1$ 
where ${\mathcal A}_1:={\mathcal A}$ and 
${\mathcal A}_i$ is the set of subsets of $S_i$ 
induced by ${\mathcal A}_{i-1}$ and the twist 
for each $i=2,\ldots, n-1$.

\begin{Definition}\label{def:compatible}
Let $(W,R)$ and $(W,S)$ be Coxeter systems and 
let ${\mathcal A}$ and ${\mathcal B}$ be 
separations of $R$ and $S$ respectively.
We say that $(W,R)$ and $(W,S)$ are 
\textit{compatible on separations ${\mathcal A}$ and ${\mathcal B}$}, 
if 
\begin{enumerate}
\item[(i)] each $A\in {\mathcal A}$ 
is conjugate to some unique $B\in {\mathcal B}$ and
\item[(ii)] each $B\in {\mathcal B}$ 
is conjugate to some unique $A\in {\mathcal A}$.
\end{enumerate}
We define that two Coxeter systems $(W,R)$ and $(W,S)$ 
are \textit{some-separation-compatible}, 
if $(W,R)$ and $(W,S)$ are 
compatible on some separations ${\mathcal A}$ and ${\mathcal B}$ of 
$R$ and $S$ respectively.
\end{Definition}

Every spherical subset $T$ of $S$ is not separated 
by any spherical(-product) subset $\sigma$ of $S$.
Hence $T$ is a twist-rigid subset of $S$ and 
$T\subset A$ for some maximal twist-rigid subset $A$ of $S$.

If two Coxeter systems are some-separation-compatible, then 
they are maximal-twist-rigid-subset-compatible and 
maximal-spherical-subset-compatible, 
hence they are angle-compatible.

We define and investigate type(I)-type(II)-compatible later.

\medskip

In this paper, we study 
when are Coxeter generating sets conjugate up to finite twists 
under ``the untangle-condition'' on conjugate subsets.

Let $(W,S)$ be a Coxeter system.
We define ``untangle-conjugate''.

\begin{Definition}\label{def:U}
For spherical subsets $\sigma$, $\tau$ and $T$ of $S$ as $\sigma \cup \tau \subset T$, 
we denote $\ds \sigma \mathop{\simeq}_{w_T} \tau$, if $w_T \sigma w_T = \tau$.

Let $U$ and $U'$ be non-empty subsets of $S$.
They are uniquely denoted by $U=U_\sigma \cup U_\nu$ and 
$U'=U'_\sigma \cup U'_\nu$.
For a spherical subset $T$ of $S$ as $U_\sigma \cup U'_\sigma \subset T$, 
we denote $\ds U \mathop{\simeq}_{w_T} U'$, 
if $w_T U_\sigma w_T = U'_\sigma$, $U_\nu = U'_\nu$ and 
$st=ts$ for any $s \in U_\nu$ and $t\in T$.

Here if $\ds U \mathop{\simeq}_{w_T} U'$, then $U$ and $U'$ are conjugate and 
\begin{align*}
w_T U w_T &= w_T (U_\sigma \cup U_\nu) w_T = (w_T U_\sigma w_T) \cup (w_T U_\nu w_T) \\
&= U'_\sigma \cup U_\nu = U'.
\end{align*}
Also if 
\[ \ds U=U_1 \mathop{\simeq}_{\,\,w_{T_1}} U_2 \mathop{\simeq}_{\,\,w_{T_2}} 
\cdots {\hspace*{-3mm}}\mathop{\simeq}_{\ \ w_{T_{q-1}}} U_q = U', \]
then $U$ and $U'$ are conjugate and for $w:=w_{T_{q-1}}\cdots w_{T_2}w_{T_1}$,
\begin{align*}
w U w^{-1} 
&= w_{T_{q-1}}\cdots w_{T_2}w_{T_1} U_1w_{T_1}w_{T_2} \cdots w_{T_{q-1}} \\
&= w_{T_{q-1}}\cdots w_{T_2} U_2 w_{T_2} \cdots w_{T_{q-1}} \\
&= \cdots \\
&= w_{T_{q-1}} U_{q-1} w_{T_{q-1}} \\
&= U_q =U'.
\end{align*}

\smallskip

We say that conjugate subsets $U$ and $U'$ of $S$ are \textit{untangle}, 
if the following statements (1) and (2) hold:
\begin{enumerate}
\item[(1)] In the case that $U \neq U'$, 
there exist a sequence $U_1,\cdots, U_q$ of subsets of $S$ and 
a sequence $T_1,\cdots, T_{q-1}$ of spherical subsets of $S$ such that 
\[ \ds U=U_1 \mathop{\simeq}_{\,\,w_{T_1}} U_2 \mathop{\simeq}_{\,\,w_{T_2}} 
\cdots {\hspace*{-3mm}}\mathop{\simeq}_{\ \  w_{T_{q-1}}} U_q = U'. \]
\item[(2)] In the case that $U=U'$, 
for $w \in W$ as $wUw^{-1}=U$, 
we define the bijective map $f_w:U \to U$ 
by $f_w(a)=w a w^{-1}$ for any $a \in U$.
For any $w \in W$ as $wUw^{-1}=U$, $f_w=f_1={\rm id}_U$ for $1\in W$, or, 
there exist a sequence $U_1,\cdots, U_q$ of subsets of $S$ and 
a sequence $T_1,\cdots, T_{q-1}$ of spherical subsets of $S$ such that 
\[ \ds U=U_1 \mathop{\simeq}_{\,\,w_{T_1}} U_2 \mathop{\simeq}_{\,\,w_{T_2}} 
\cdots {\hspace*{-3mm}}\mathop{\simeq}_{\ \ w_{T_{q-1}}} U_q = U \]
and $f_w = f_{w_0}$ for $w_0:=w_{T_{q-1}}\cdots w_{T_2}w_{T_1}$.
\end{enumerate}
\end{Definition}

Now we define ``the untangle-conjugate-condition'' and 
``the untangle-condition''.

\begin{Definition}\label{def:untangle-condition}
We say that a Coxeter system $(W,S)$ has 
\textit{the untangle-conjugate-condition}, if 
all conjugate subsets $U$ and $U'$ of $S$ are untangle 
and the following holds:
\begin{enumerate}
\item[($*$)] Let $R$ be a Coxeter generating set for $W$, 
let $S_1,\ldots,S_n$ be the connected components of $S$ 
and let $R_1,\ldots,R_n$ be the connected components of $R$ 
(that is, $W=W_{S_1}*\cdots *W_{S_n}=W_{R_1}*\cdots *W_{R_n}$).
If $S_i$ and $R_i$ are conjugate in $W$ for any $i=1,\ldots,n$, 
then there exists a Coxeter generating set $R'$ for $W$ 
such that $R'$ is obtained from $R$ by some finite twists and $S$ and $R'$ are conjugate.
\end{enumerate}
Also we say that 
a Coxeter system $(W,S)$ has \textit{the untangle-condition}, if 
all Coxeter systems obtained from $(W,S)$ by some finite twists 
have the untangle-conjugate-condition.
\end{Definition}

\begin{Example}\label{example00}
We consider a Coxeter system $(W,S)$ defined by Figure~\ref{figA}~(i) or (ii).
\begin{figure}[H]
\centering
{ 
\medskip
\includegraphics[keepaspectratio, scale=0.90, bb=0 0 317 77]{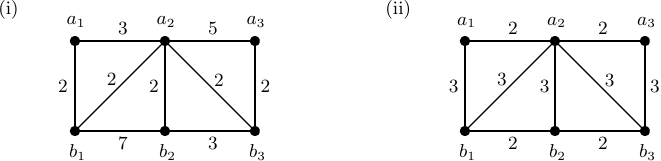} \ 
}
\caption{Example~\ref{example00}}\label{figA}
\end{figure}

Here $S$ is the vertex set of the figure and $o(ss)=1$ for any $s\in S$.
If two vertices $s,t \in S$ as $s\neq t$ 
do not span any edge in the figure then we define $o(st)=\infty$.
If two vertices $s,t \in S$ as $s\neq t$ 
span an edge numbering $m$ then define $o(st)=m$.

Let $\sigma=\sigma_1:=\{ a_1,b_1 \}$, $\sigma_2:=\{ a_2,b_1 \}$, $\sigma_3:=\{ a_2,b_2 \}$, 
$\sigma_4:=\{ a_2,b_3 \}$ and $\tau=\sigma_5:=\{ a_3,b_3 \}$.
Also let $T_1:=\{ a_1,a_2,b_1 \}$, $T_2:=\{ a_2,b_1,b_2 \}$, 
$T_3:=\{ a_2,b_2,b_3 \}$ and $T_4:=\{ a_2,a_3,b_3 \}$.

Then $\sigma$ and $\tau$ are untangle-conjugate and 
\[ \ds \sigma = \sigma_1 \mathop{\simeq}_{\,\,w_{T_1}} \sigma_2 \mathop{\simeq}_{\,\,w_{T_2}} 
\sigma_3 \mathop{\simeq}_{\,\,w_{T_3}} \sigma_4 \mathop{\simeq}_{\,\,w_{T_4}} \sigma_5= \tau \]
in both cases (i) and (ii).

\smallskip

In the case (ii), for $T_0:=\sigma =\{ a_1,b_1 \}$, 
\[ \ds \sigma  \mathop{\simeq}_{\,\,w_{T_0}} \sigma. \]
Then $w_{T_0} \sigma w_{T_0} =\sigma$. Here $w_{T_0} a_1 w_{T_0} = b_1$ and $w_{T_0} b_1 w_{T_0} = a_1$.
Hence $f_{w_{T_0}}:\sigma \to \sigma$ is the bijective map such that 
$f_{w_{T_0}}(a_1)=b_1$ and $f_{w_{T_0}}(b_1)=a_1$.
\end{Example}

We can see \cite{BH} and \cite{De} on conjugate spherical subsets.

\begin{Example}\label{example001}
We consider a Coxeter system $(W,S)$ defined by Figure~\ref{figB}.
\begin{figure}[H]
\centering
{ 
\medskip
\hspace*{5mm}\includegraphics[keepaspectratio, scale=0.90, , bb=0 0 207 70]{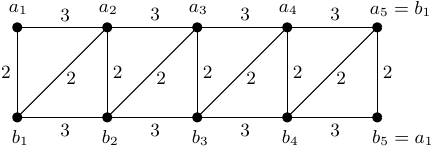}
}
\caption{Example~\ref{example001}}\label{figB}
\end{figure}

Here we identify $a_1=b_5$ and $b_1=a_5$.
(The nerve of $(W,S)$ is a triangulation of the M\"{o}bius band.)

Let $\sigma=\sigma_1:=\{ a_1,b_1 \}=\{ b_5,a_5 \}$, $\sigma_2:=\{ a_2,b_1 \}$, $\sigma_3:=\{ a_2,b_2 \}$, 
$\sigma_4:=\{ a_3,b_2 \}$, $\sigma_5:=\{ a_3,b_3 \}$, 
$\sigma_6:=\{ a_4,b_3 \}$, $\sigma_7:=\{ a_4,b_4 \}$ and $\sigma_8:=\{ a_5,b_4 \}$.
Also let $T_1:=\{ a_1,a_2,b_1 \}$, $T_2:=\{ a_2,b_1,b_2 \}$, 
$T_3:=\{ a_2,a_3,b_2 \}$, $T_4:=\{ a_3,b_2,b_3 \}$, 
$T_5:=\{ a_3,a_4,b_3 \}$, $T_6:=\{ a_4,b_3,b_4 \}$, 
$T_7:=\{ a_4,a_5,b_4 \}$ and $T_8:=\{ a_5,b_4,b_5 \}$.

Then 
\[ \ds \sigma = \sigma_1 \mathop{\simeq}_{\,\,w_{T_1}} \sigma_2 \mathop{\simeq}_{\,\,w_{T_2}} 
\sigma_3 \mathop{\simeq}_{\,\,w_{T_3}} \sigma_4 \mathop{\simeq}_{\,\,w_{T_4}} 
\sigma_5 \mathop{\simeq}_{\,\,w_{T_5}} \sigma_6 \mathop{\simeq}_{\,\,w_{T_6}} 
\sigma_7 \mathop{\simeq}_{\,\,w_{T_7}} \sigma_8 \mathop{\simeq}_{\,\,w_{T_8}} \sigma_1 =\sigma. \]
For $w_0:=w_{T_8}w_{T_7}w_{T_6}w_{T_5}w_{T_4}w_{T_3}w_{T_2}w_{T_1}$, we have that $w_0 \sigma w_0^{-1}=\sigma$.
Here $w_0 a_1 w_0^{-1}=b_1$ and $w_0 b_1 w_0^{-1}=a_1$.

For example, 
if $w=1$, $w=a_1$, $w=b_1$, $w=w_{\sigma}=a_1b_1$, $w=w_0$, $w=w_0^2$, or $w=a_1 w_0 a_1 b_1 w_0$, 
then $w \sigma w^{-1}=\sigma$.

For any $w\in W$ as $w \sigma w^{-1}=\sigma$, 
if [\,$w a_1 w^{-1}=a_1$ and $w b_1 w^{-1} =b_1$\,] then $f_w = f_1$ for $1\in W$.
Also if [\,$w a_1 w^{-1}=b_1$ and $w b_1 w^{-1} =a_1$\,] then $f_w = f_{w_0}$ as above.
\end{Example}

The author does not know whether there is an example of 
non-untangle-conjugate subsets of a Coxeter generating set.

After some preliminaries in Sections~\ref{sec2} and \ref{sec3}, 
we prove the following theorem in Section~\ref{sec4}.

\begin{Theorem}\label{MainTheorem}
Let $(W,R)$ and $(W,S)$ be Coxeter systems with the untangle-condition.
If $(W,R)$ and $(W,S)$ are some-separation-compatible, then 
$R$ and $S$ are conjugate up to finite twists.
\end{Theorem}

We obtain a corollary from Theorem~\ref{MainTheorem}.

\begin{Corollary}\label{Cor}
For Coxeter systems $(W,R)$ and $(W,S)$ with the untangle-condition, 
the following statements are equivalent$:$
\begin{enumerate}
\item[\textnormal{(i)}] $R$ and $S$ are conjugate up to finite twists.
\item[\textnormal{(ii)}] $(W,R_0)$ and $(W,S)$ are some-separation-compatible 
for some Coxeter generating set $R_0$ obtained from $R$ by finite twists.
\end{enumerate}
\end{Corollary}

Now we define ``type(I)'' and ``type(II)'' subsets 
and ``the standard-separation'' of a Coxeter generating set.

\begin{Definition}\label{def:standard-separation}
Let $(W,S)$ be a Coxeter system.
For each minimal separation ${\mathcal A}$ of $S$, 
we define ${\mathcal U}_{\mathcal A}$ 
as the set of separators of ${\mathcal A}$.
Let 
\[ \overline{\mathcal U}:=
\bigcap\{{\mathcal U}_{\mathcal A}:
\text{${\mathcal A}$ is a minimal separation of $S$} \}. \]
Then we define the separation $\widetilde{\mathcal A}_S$ of $S$ 
by $\overline{\mathcal U}$ as follows.

\smallskip

Let ${\mathcal A}_0$ be the set of maximal twist-rigid subsets of $S$.
For $A,A'\in {\mathcal A}_0$, 
we denote $A \sim A'$ if any $U\in \overline{\mathcal U}$ does not separate 
$A$ and $A'$; that is, 
there exist a minimal separation ${\mathcal A}$ of $S$ and 
an element $\overline{A}\in {\mathcal A}$ such that $A\cup A' \subset \overline{A}$.
Then ``$\sim$'' is an equivalence relation on the set ${\mathcal A}_0$.
Let $[A]:=\{A'\in {\mathcal A}_0:A'\sim A\}$ 
that is the equivalence class for $A\in {\mathcal A}_0$.
Here 
\[ {\mathcal A}_0/\!\! \sim \; =\{[A]: A\in {\mathcal A}_0 \}=\{ [A_1],\ldots,[A_n] \} \]
for some $A_1,\ldots,A_n\in {\mathcal A}_0$ as $[A_i] \neq [A_j]$ if $i\neq j$.
Let 
\[ \overline{A}_i:=\bigcup [A_i]=\bigcup \{A\in {\mathcal A}_0 : A\sim A_i \} \]
for each $i=1,\ldots,n$.
Then we define 
\[ \widetilde{\mathcal A}_S:=\{ \overline{A}_1,\ldots,\overline{A}_n \}. \]
In Section~\ref{sec5}, we show that $\widetilde{\mathcal A}_S$ is a separation of $S$.
We say that $\widetilde{\mathcal A}_S$ is the \textit{standard separation} of $S$.

A subset $A$ of $S$ is said to be \textit{type(I)}, if 
$A\in {\mathcal A}$ holds for any minimal separation ${\mathcal A}$ of $S$.
Let $\widetilde{\mathcal A}_S^{\rm (I)}$ be the set of type(I) subsets of $S$; that is,
\[ \widetilde{\mathcal A}_S^{\rm (I)}=\{ A\in \widetilde{\mathcal A}_S : 
A\in {\mathcal A} \ \text{for any minimal separation ${\mathcal A}$ of $S$} \}. \]
We also define 
\[ \widetilde{\mathcal A}_S^{\rm (II)}:= \widetilde{\mathcal A}_S - \widetilde{\mathcal A}_S^{\rm (I)} \]
and every element $A\in \widetilde{\mathcal A}_S^{\rm (II)}$ is called a \textit{type(II)} subset of $S$.

Then $\widetilde{\mathcal A}_S = \widetilde{\mathcal A}_S^{\rm (I)} \cup \widetilde{\mathcal A}_S^{\rm (II)}$ 
is a disjoint union.
Here $\widetilde{\mathcal A}_S$, $\widetilde{\mathcal A}_S^{\rm (I)}$ and $\widetilde{\mathcal A}_S^{\rm (II)}$ 
are uniquely determined by the Coxeter system $(W,S)$.
\end{Definition}

Let $(W,S)$ be a Coxeter system.
We say that for a subset $A$ of $S$, 
$A'$ is a twist of $A$ that \textit{induces} some twist of $S$, 
if there exist a spherical-product subset $U$ of $S$ and $w\in W$ 
such that $U\subset A$, $U$ separates $A$, 
$A'$ is a twist of $A$ obtained by $U$ and $w$, 
$U$ separates $S$ and some twist $S'$ of $S$ is obtained by $U$ and $w$.
We also say that 
$A''$ is obtained from $A$ by some finite twists 
that \textit{induces} some twist of $S$, 
if there exist a sequence $A=A_1,A_2,\ldots,A_n=A''$ and 
a sequence $S=S_1,S_2,\ldots,S_n=S''$ of Coxeter generating sets for $W$ 
such that $A_i \subset S_i$ for any $i=1,\ldots, n$ and 
$A_{i+1}$ is a twist of $A_i$ 
that induces a twist $S_{i+1}$ of $S_i$ for each $i=1,\ldots, n-1$.

We define ``type(II)-compatible''.

\begin{Definition}\label{def:type(II)-compatible}
Let $(W,R)$ and $(W,S)$ be Coxeter systems.
We say that $A\in \widetilde{\mathcal A}_R^{\rm (II)}$ and 
$B\in \widetilde{\mathcal A}_S^{\rm (II)}$ are \textit{type(II)-compatible}, 
if 
\begin{enumerate}
\item[~] 
there exists $A_0$ obtained from $A$ by some finite twists 
that induce some twists of $R$ preserving $\widetilde{\mathcal A}_R-\{ A \}$ 
such that $A_0$ and $B$ are conjugate in $W$.
\end{enumerate}
This means that 
the type(II) subsets $A\in \widetilde{\mathcal A}_R^{\rm (II)}$ and 
$B\in \widetilde{\mathcal A}_S^{\rm (II)}$ 
are conjugate up to finite twists that induce some twists of $R$ 
preserving $\widetilde{\mathcal A}_R-\{ A \}$.
\end{Definition}

Now we define ``type(I)-type(II)-compatible''.

\begin{Definition}\label{def:type(I)-type(II)-compatible}
Two Coxeter systems $(W,R)$ and $(W,S)$ 
are said to be \textit{type(I)-type(II)-compatible}, 
if for the standard separations $\widetilde{\mathcal A}_R$ and $\widetilde{\mathcal A}_S$ 
of $R$ and $S$ respectively, 
\begin{enumerate}
\item[(i)] each $A\in \widetilde{\mathcal A}_R^{\rm (I)}$ 
is conjugate to some unique $B\in \widetilde{\mathcal A}_S^{\rm (I)}$, 
\item[(ii)] each $B\in \widetilde{\mathcal A}_S^{\rm (I)}$ 
is conjugate to some unique $A\in \widetilde{\mathcal A}_R^{\rm (I)}$, 
\item[(iii)] for each $A\in \widetilde{\mathcal A}_R^{\rm (II)}$, 
there exists a unique $B\in \widetilde{\mathcal A}_S^{\rm (II)}$ 
such that $A$ and $B$ are type(II)-compatible, and
\item[(iv)] for each $B\in \widetilde{\mathcal A}_S^{\rm (II)}$, 
there exists a unique $A\in \widetilde{\mathcal A}_R^{\rm (II)}$ 
such that $B$ and $A$ are type(II)-compatible.
\end{enumerate}
\end{Definition}

In Section~\ref{sec6}, we show that 
if two Coxeter systems with the untangle-condition are type(I)-type(II)-compatible, then 
they are some-separation-compatible up to finite twists.
Thus, if two Coxeter systems $(W,R)$ and $(W,S)$ with the untangle-condition are type(I)-type(II)-compatible, 
then $R$ and $S$ are conjugate up to finite twists by Theorem~\ref{MainTheorem}.

We show the following theorem in Section~\ref{sec6}.

\begin{Theorem}\label{MainTheorem2}
For Coxeter systems $(W,R)$ and $(W,S)$ with the untangle-condition, 
the following two statements are equivalent$:$
\begin{enumerate}
\item[\textnormal{(i)}] $R$ and $S$ are conjugate up to finite twists.
\item[\textnormal{(ii)}] $(W,R)$ and $(W,S)$ are type(I)-type(II)-compatible.
\end{enumerate}
\end{Theorem}

For given Coxeter systems $(W,R)$ and $(W,S)$, 
it seems that to consider whether
\begin{enumerate}
\item[(a)] $(W,R)$ and $(W,S)$ are type(I)-type(II)-compatible
\end{enumerate}
is more simple than to consider whether
\begin{enumerate}
\item[(b)] $(W,R_0)$ and $(W,S)$ are some-separation-compatible 
for some Coxeter generating set $R_0$ obtained from $R$ by finite twists.
\end{enumerate}

In Theorem~\ref{MainTheorem2}, we use two conditions as ``the untangle-condition'' and 
``type(I)-type(II)-compatible'' for Coxeter systems.

\begin{Problem}\label{ProblemI}
The untangle-conjugate-condition will always hold for all Coxeter systems $(W,S)$?
\end{Problem}

\begin{Problem}\label{ProblemII}
Angle-compatible Coxeter systems $(W,R)$ and $(W,S)$ will be always type(I)-type(II)-compatible?
\end{Problem}

If there exist counter-examples of 
angle-compatible Coxeter systems $(W,R)$ and $(W,S)$ 
such that any $R_0$ and $S_0$ obtained from $R$ and $S$ by some finite twists respectively 
are not type(I)-type(II)-compatible, or 
a Coxeter system $(W,S)$ that does not have the untangle-conjugate-condition, 
then they are meaningful examples.

If Problems~\ref{ProblemI} and \ref{ProblemII} both can be solved affirmatively, 
then the twist-conjecture for Coxeter systems can be solved affirmatively from Theorem~\ref{MainTheorem2}.
Also if Problem~\ref{problem1.2} can be solved from the twist-conjecture for Coxeter systems affirmatively, 
then the isomorphism problem for Coxeter groups of finite ranks can be solved by Theorem~\ref{thmMM}.

\section{Remarks and examples on separations}\label{sec2}

We introduce some remarks and examples on 
separations of Coxeter generating sets.

\begin{Remark}\label{rem3-0}
Let $(W,S)$ be a Coxeter system as $S$ is connected.
Let ${\mathcal A}$ be a separation of $S$ and let $n:=|{\mathcal A}|$.
Then we may denote ${\mathcal A}=\{A_1,A_2,\ldots,A_n\}$ 
such that for each $i=1,\ldots,n-1$, 
$U_i:=(A_1\cup \cdots \cup A_i) \cap A_{i+1}$ is 
maximal in 
\[ \{(A_1\cup \cdots \cup A_i) \cap A_j : j=i+1,\ldots,n \}. \]
Here by the definition of a separation of $S$, 
$U_i$ is a separator of ${\mathcal A}$ and 
it is a spherical-product subset of $S$ 
that separates $S$ for any $i=1,\ldots,n-1$.
Hence $W$ has a structure as 
\[ W=(\cdots ((W_{A_1} *_{W_{U_1}} W_{A_2})*_{W_{U_2}} W_{A_3})*_{W_{U_3}}\cdots 
)*_{W_{U_{n-1}}} W_{A_n}. \]
\end{Remark}

\begin{Remark}\label{Rem3-2}
Let $(W,S)$ be a Coxeter system as $S$ is connected and 
let $U$ be a spherical-product subset of $S$ that separates $S$.
Let ${\mathcal A}_0$ be the set of maximal twist-rigid subsets of $S$.
Suppose that $X_1, \ldots, X_t$ are the connected components of $S-U$ and 
$S-U= X_1 \cup \cdots \cup X_t$ is a disjoint union.
Let 
\begin{align*}
&\overline{X}_i:=\bigcup \{ A\in {\mathcal A}_0 : A \subset X_i\cup U 
\ \text{and} \ A\not\subset U \} \\
\intertext{for each $i=1,\ldots,t$ and let }
&\overline{Y}:=\bigcup \{ A\in {\mathcal A}_0 : A \subset U \}. 
\end{align*}
Here each $\overline{X}_i$ (and $\overline{Y}$) is connected and it is 
a union of some maximal twist-rigid subsets of $S$ (if $\overline{Y}$ is non-empty).
Let ${\mathcal A}:=\{\overline{X}_1,\ldots , \overline{X}_t, \overline{Y} \}$ 
if $\overline{Y}$ is non-empty, and let 
${\mathcal A}:=\{\overline{X}_1,\ldots , \overline{X}_t \}$ 
if $\overline{Y}$ is empty.
Then ${\mathcal A}$ is a separation of $S$.
We say that ${\mathcal A}$ is the induced separation of $S$ by $U$.

Let ${\mathcal A}'$ be a minimal separation of $S$ 
such that ${\mathcal A}' \preceq {\mathcal A}$.
Then $U$ does not separate any $A\in {\mathcal A}'$.
Thus we can obtain a minimal separation ${\mathcal A}'$ 
of $S$ from a spherical-product subset $U$ of $S$ 
that separates $S$ such that $U$ does not separate any $A\in {\mathcal A}'$.

\smallskip

Let $\{U_j\}$ be the set of spherical-product subsets of $S$ that separate $S$.
By the above argument, 
we obtain the separation ${\mathcal A}_j$ of $S$ induced by each $U_j$.
If $\{U_j\}$ is non-empty 
(that is, if $(W,S)$ is not a twist-rigid Coxeter system), 
then for any minimal separation ${\mathcal A}'$ of $S$,
${\mathcal A}' \preceq {\mathcal A}_j$ for some $j$.
Hence, from considering the set $\{U_j\}$ 
of spherical-product subsets that separate $S$ 
and the induced separations ${\mathcal A}_j$, 
we can obtain all minimal separations of $S$.
\end{Remark}

We give some examples.

\begin{Example}\label{example1B}
We consider Coxeter systems 
$(W_i,S_i)$ ($i=1,2,3,4$) defined by Figure~\ref{fig1B}.

\begin{figure}[h]
\centering
{ 
\includegraphics[keepaspectratio, scale=0.95, bb=0 0 351 120]{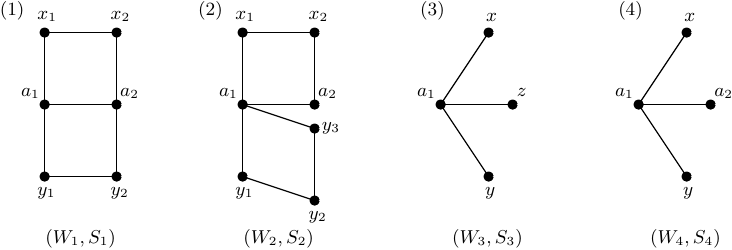}
}
\caption{Example~\ref{example1B}}\label{fig1B}
\end{figure}

Here $S_1$, $S_2$, $S_3$ and $S_4$ are the vertex sets of the figures 
$(1)$, $(2)$, $(3)$ and $(4)$ respectively in Figure~\ref{fig1B}.
Let $i\in \{1,2,3,4 \}$.
We define $o(ss)=1$ for any $s\in S_i$.
If two vertices $s,t \in S_i$ as $s\neq t$ 
do not span any edge in the figure 
then we define $o(st)=\infty$.
If two vertices $s,t \in S_i$ as $s\neq t$ 
span an edge then we consider that 
$o(st)=o(ts)$ has some (arbitrary) finite number at least $2$.
Then a Coxeter system $(W_i,S_i)$ is obtained.

(1) Let $U:=\{a_1,a_2\}$, $X_1:=\{ x_1,x_2\}$ and $X_2:=\{ y_1,y_2\}$ in $(W_1,S_1)$.
Then $U$ is a spherical-product subset of $S_1$ that separates $S_1$ and 
$S_1-U=X_1 \cup X_2$ where $X_1$ and $X_2$ are the connected components of $S_1-U$.
Then $\overline{X}_1=\{ x_1,x_2, a_1,a_2\}$ and $\overline{X}_2=\{ y_1,y_2, a_1,a_2\}$.
Here ${\mathcal A}:=\{\overline{X}_1,\overline{X}_2 \}$ is the separation of $S_1$ 
induced by $U$, and $U$ is a separator of ${\mathcal A}$.

(2) Let $U:=\{a_1,a_2\}$, $X_1:=\{ x_1,x_2\}$ and $X_2:=\{ y_1,y_2,y_3\}$ in $(W_2,S_2)$.
Then $U$ is a spherical-product subset of $S_2$ that separates $S_2$ and 
$S_2-U=X_1 \cup X_2$ where $X_1$ and $X_2$ are the connected components of $S_2-U$.
Then $\overline{X}_1=\{ x_1,x_2, a_1,a_2\}$ and $\overline{X}_2=\{ y_1,y_2,y_3, a_1 \}$.
Here ${\mathcal A}:=\{\overline{X}_1,\overline{X}_2 \}$ is the separation of $S_2$ 
induced by $U$, and $U':=\{ a_1 \}$ is a separator of ${\mathcal A}$ 
(here $U=\{ a_1,a_2 \}$ is not a separator of ${\mathcal A}$).

(3) Let $U:=\{a_1\}$, $X_1:=\{ x\}$, $X_2:=\{ y\}$ and $X_3:=\{ z\}$ in $(W_3,S_3)$.
Then $U$ is a spherical-product subset of $S_3$ that separates $S_3$ and 
$S_3-U=X_1 \cup X_2\cup X_3$ where $X_1$, $X_2$ and $X_3$ are the connected components of $S_3-U$.
Then $\overline{X}_1=\{ x, a_1\}$, $\overline{X}_2=\{ y, a_1\}$ and $\overline{X}_3=\{ z, a_1 \}$.
Here ${\mathcal A}:=\{\overline{X}_1,\overline{X}_2,\overline{X}_3 \}$ is the separation of $S_3$ 
induced by $U$, and $U$ is a separator of ${\mathcal A}$.

(4) Let $U:=\{a_1,a_2\}$, $X_1:=\{ x\}$ and $X_2:=\{ y\}$ in $(W_4,S_4)$.
Then $U$ is a spherical-product subset of $S_4$ that separates $S_4$ and 
$S_4-U=X_1 \cup X_2$ where $X_1$ and $X_2$ are the connected components of $S_4-U$.
Then $\overline{X}_1=\{ x, a_1\}$, $\overline{X}_2=\{ y, a_1\}$ and $\overline{Y}=\{ a_1,a_2 \}$.
Here ${\mathcal A}:=\{\overline{X}_1,\overline{X}_2,\overline{Y} \}$ is the separation of $S_4$ 
induced by $U$, and $U':=\{a_1\}$ is a separator of ${\mathcal A}$ 
(here $U=\{ a_1,a_2 \}$ is not a separator of ${\mathcal A}$).
\end{Example}

\begin{Example}\label{example11}
We consider a Coxeter system $(W,S)$ defined by Figure~\ref{fig11}.
\begin{figure}[H]
\centering
{ 
\medskip
\includegraphics[keepaspectratio, scale=0.90, bb=0 0 188 188]{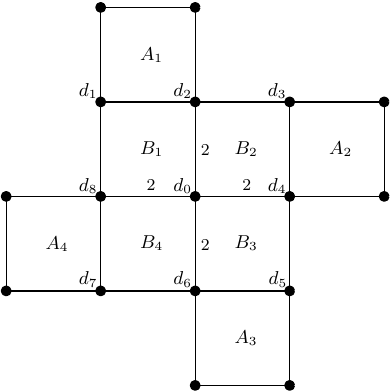}
}
\caption{Example~\ref{example11}}\label{fig11}
\end{figure}

Here $S$ is the vertex set of the figure.
Let $o(ss)=1$ for any $s\in S$.
If two vertices $s,t \in S$ as $s\neq t$ 
do not span any edge in the figure 
then we define $o(st)=\infty$.
If two vertices $s,t \in S$ as $s\neq t$ 
span an edge numbering $2$ then define $o(st)=2$.
Also if two vertices $s,t \in S$ as $s\neq t$ 
span an edge with no numbering then we consider that 
$o(st)=o(ts)$ has some (arbitrary) finite number at least $2$.
Then a Coxeter system $(W,S)$ is obtained.

Let $A_i$ be the set of $4$-vertices around ``$A_i$'' in the figure 
and let $B_i$ be the set of $4$-vertices around ``$B_i$'' 
in the figure for each $i=1,2,3,4$.
\begin{enumerate}
\item[(a)] The set of maximal twist-rigid subsets of $S$ is 
\[ \{A_1,\ A_2,\ A_3,\ A_4,\  B_1,\ B_2,\ B_3,\ B_4\}. \]
Here 
\[ D:=\{d_1,\,d_2,\,d_3,\,d_4,\,d_5,\,d_6,\,d_7,\,d_8 \} \]
is not a twist-rigid subset of $S$.
Indeed for example $U:=\{ d_2,\, d_0,\, d_6 \}$ 
is a spherical-product subset of $S$ that separates $S$ and 
$U\cap D=\{d_2,\, d_6 \}$ separates $D$.
We also note that $(W_D,D)$ is a twist-rigid Coxeter system.
\item[(b)] The minimal separations of $S$ are the 4-sets as 
\begin{enumerate}
\item[(1)] $\{A_1,\ A_2,\ A_3,\ A_4,\  B_1 \cup B_2 ,\  B_3\cup B_4 \}$,
\item[(2)] $\{A_1,\ A_2,\ A_3,\ A_4,\  B_1 \cup B_4 ,\  B_2\cup B_3 \}$,
\item[(3)] $\{A_1,\ A_2,\ A_3,\ A_4,\  B_1\cup B_3,\  B_2,\  B_4 \}$ and
\item[(4)] $\{A_1,\ A_2,\ A_3,\ A_4,\  B_1,\  B_3,\  B_2\cup B_4 \}$.
\end{enumerate}
\item[(c)] The set of type(I) subsets of $S$ is 
\[ \widetilde{\mathcal A}_S^{\rm (I)}=\{A_1,\ A_2,\ A_3,\ A_4\}. \]
\item[(d)] The set of type(II) subsets of $S$ is 
\[ \widetilde{\mathcal A}_S^{\rm (II)}=\{B_1\cup B_2\cup B_3\cup B_4\}. \]
\item[(e)] The standard separation of $S$ is 
\begin{align*}
\widetilde{\mathcal A}_S &= \widetilde{\mathcal A}_S^{\rm (I)}\cup \widetilde{\mathcal A}_S^{\rm (II)} \\
&= \{A_1,\ A_2,\ A_3,\ A_4,\  B_1\cup B_2\cup B_3\cup B_4\}. 
\end{align*}
\end{enumerate}
\end{Example}

\begin{Example}\label{example13}
We consider a Coxeter system $(W,S)$ defined by Figure~\ref{fig13}.
\begin{figure}[H]
\centering
{ 
\medskip
\includegraphics[keepaspectratio, scale=0.90, bb=0 0 192 99]{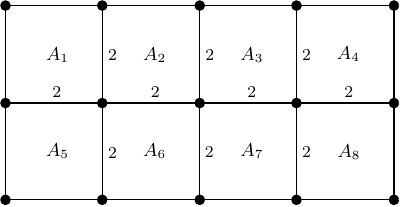}
}
\caption{Example~\ref{example13}}\label{fig13}
\end{figure}

Here $S$ is the vertex set of the figure and 
$o(st)$ is defined as in Example~\ref{example11}.

Let $A_i$ be the set of $4$-vertices around ``$A_i$'' 
in the figure for each $i=1,2,3,4,5,6,7,8$.
\begin{enumerate}
\item[(a)] The set of maximal twist-rigid subsets of $S$ is 
\[ \{A_1,\ A_4,\ A_5,\ A_8,\ A_2\cup A_6,\ A_3\cup A_7 \}. \]
\item[(b)] The minimal separations of $S$ are the 9-sets as 
\begin{enumerate}
\item[(1)] $\{A_1\cup A_5,\ A_2\cup A_6,\  A_3\cup A_7,\ A_4\cup A_8 \}$,
\item[(2)] $\{A_1\cup A_5,\ A_2\cup A_6,\   A_3\cup A_7\cup A_8,\ A_4 \}$,
\item[(3)] $\{A_1\cup A_5,\ A_2\cup A_6,\    A_3\cup A_4\cup A_7,\ A_8 \}$,
\item[(4)] $\{A_1,\ A_2\cup A_5\cup A_6,\  A_3\cup A_7,\ A_4\cup A_8 \}$,
\item[(5)] $\{A_1,\ A_2\cup A_5\cup A_6,\   A_3\cup A_7\cup A_8,\ A_4 \}$,
\item[(6)] $\{A_1,\ A_2\cup A_5\cup A_6,\    A_3\cup A_4\cup A_7,\ A_8 \}$,
\item[(7)] $\{A_1\cup A_2\cup A_6,\ A_5,\  A_3\cup A_7,\ A_4\cup A_8 \}$,
\item[(8)] $\{A_1\cup A_2\cup A_6,\ A_5,\   A_3\cup A_7\cup A_8,\ A_4 \}$ and
\item[(9)] $\{A_1\cup A_2\cup A_6,\ A_5,\    A_3\cup A_4\cup A_7,\ A_8 \}$.
\end{enumerate}
\item[(c)] The set $\widetilde{\mathcal A}_S^{\rm (I)}$ of type(I) subsets of $S$ is empty.
\item[(d)] The set of type(II) subsets of $S$ is 
\[ \widetilde{\mathcal A}_S^{\rm (II)}=\{ B_1,\ B_2 \} \]
where $B_1:=A_1\cup A_2 \cup A_5\cup A_6$ and 
$B_2:=A_3\cup A_4 \cup A_7\cup A_8$.
\item[(e)] The standard separation of $S$ is 
\[ \widetilde{\mathcal A}_S = \{ B_1,\ B_2 \}. \]
\end{enumerate}
\end{Example}

\begin{Example}\label{example14}
We consider the Coxeter system $(W,S)$ defined by Figure~\ref{fig14}.
Here $S:=\{a,\,b_1,\,b_2,\,b_3,\,c\}$.
\begin{figure}[H]
\centering
{ 
\includegraphics[keepaspectratio, scale=0.90, bb=0 0 131 119]{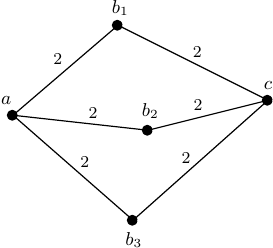}
}
\caption{Example~\ref{example14}}\label{fig14}
\end{figure}

Let $A_1:= \{ a,\,b_1,\,c \}$, $A_2:= \{ a,\,b_2,\,c \}$, 
$A_3:= \{ a,\,b_3,\,c \}$, $B_1:= \{ a,\,b_1,\,b_2,\,b_3 \}$ and 
$B_2:= \{ b_1,\,b_2,\,b_3,\,c \}$.
\begin{enumerate}
\item[(a)] The set of maximal twist-rigid subsets of $S$ is 
\[ \{\, \{ a,\,b_1 \},\, \{ a,\,b_2 \},\,\{ a,\,b_3 \},\, 
\{ b_1,\,c \},\, \{ b_2,\,c \},\,\{ b_3,\,c \}\, \}. \]
\item[(b)] The minimal separations of $S$ are the 2-sets as 
\begin{enumerate}
\item[(1)] $\{A_1,\ A_2,\ A_3 \}$ and 
\item[(2)] $\{B_1,\ B_2 \}$.
\end{enumerate}
\item[(c)] The set $\widetilde{\mathcal A}_S^{\rm (I)}$ of type(I) subsets of $S$ is empty.
\item[(d)] The set of type(II) subsets of $S$ is $\widetilde{\mathcal A}_S^{\rm (II)}=\{ S \}$.
\item[(e)] The standard separation of $S$ is $\widetilde{\mathcal A}_S = \{ S \}$.
\end{enumerate}
\end{Example}

\begin{Example}\label{example15}
We consider a Coxeter system $(W,S)$ defined by Figure~\ref{fig15}.
Here $S$ is the vertex set of the figure; that is, 
\[ S:=\{a,\,b,\,c,\,d,\,e,\,f,\,g,\,h,\,i,\,j,\,k,\,l \}, \]
and we define $o(ek)\ge 3$.

\begin{figure}[H]
\centering
{ 
\includegraphics[keepaspectratio, scale=0.90, bb=0 0 204 115]{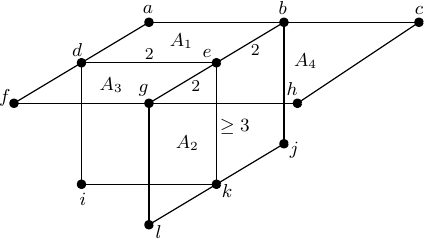}
}
\caption{Example~\ref{example15}}\label{fig15}
\end{figure}

Let 
\begin{align*}
A_1&:= \{ a,\,b,\,d,\,e \}, \\
A_2&:= \{ b,\,d,\,e,\,g,\,i,\,j,\,k,\,l \}, \\
A_3&:= \{ d,\,e,\,f,\,g \} \ \text{and} \\
A_4&:= \{ b,\,c,\,e,\,g,\,h \}.
\end{align*}
Then 
\[ {\mathcal A}:=\{ A_1,\ A_2,\ A_3,\ A_4 \} \]
is a separation of $S$.
For example, 
\begin{align*}
&U_1:=A_1\cap A_2 \ \text{is maximal in}\ \{ A_1\cap A_j : j=2,3,4 \}, \\
&U_2:=(A_1\cup A_2)\cap A_3 \ \text{is maximal in}\ \{ (A_1\cup A_2)\cap A_j : j=3,4 \} \ \text{and} \\
&U_3:=(A_1\cup A_2\cup A_3)\cap A_4 \ \text{is maximal in}\ \{ (A_1\cup A_2\cup A_3)\cap A_j : j=4 \}.
\end{align*}
Also $U_1$, $U_2$ and $U_3$ are separators of ${\mathcal A}$ and they are 
spherical-product subsets of $S$ that separate $S$.
Here we note that $A_1\cap A_3$, $A_1\cap A_4$ and $A_3\cap A_4$ are 
spherical-product subsets that do not separate $S$.

The spherical-product subsets of $S$ that separate $S$ are 
the 3-sets as 
\[U_1=\{ b,\, d,\, e \},\ 
U_2=\{ d,\, e,\, g \} \ \text{and}\ 
U_3=\{ b,\, e,\, g \}. \]
Each $U_i$ induces the separation ${\mathcal A}_i$ of $S$ 
as in Remark~\ref{Rem3-2}.
Then for any minimal separation ${\mathcal A}'$ of $S$, 
${\mathcal A}' \preceq {\mathcal A}_i$ for some $i=1,2,3$.
Here ${\mathcal A} \preceq {\mathcal A}_i$ holds for any $i=1,2,3$, 
where ${\mathcal A}=\{ A_1,A_2,A_3,A_4 \}$.
Hence ${\mathcal A}$ is the unique minimal separation of $S$.
Thus, there are no type(II) subsets of $S$ and 
\[ \widetilde{\mathcal A}_{S}= \widetilde{\mathcal A}_{S}^{\rm (I)}=\{A_1,\ A_2,\ A_3,\ A_4 \}. \]
\end{Example}

\begin{Example}\label{example16}
We consider a Coxeter system $(W',S')$ defined by Figure~\ref{fig16}.
Here $S'$ is the vertex set of the figure; that is, 
\[ S':=\{a,\,b,\,c,\,d,\,e,\,f,\,g,\,h \}. \]

\begin{figure}[H]
\centering
{ 
\includegraphics[keepaspectratio, scale=0.90, bb=0 0 205 52]{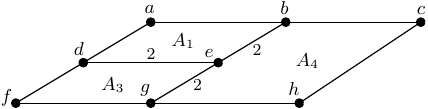}
}
\caption{Example~\ref{example16}}\label{fig16}
\end{figure}

Then $W'$ can be considered 
as a standard subgroup of $(W,S)$ in Example~\ref{example15} 
generated by $S'$.

Let 
\begin{align*}
A_1&:= \{ a,\,b,\,d,\,e \}, \\
A_3&:= \{ d,\,e,\,f,\,g \} \ \text{and} \\
A_4&:= \{ b,\,c,\,e,\,g,\,h \}.
\end{align*}
Then the spherical-product subsets of $S'$ that separate $S'$ 
are the 3-sets as 
\[U_1=\{ b,\, d,\, e \},\ 
U_2=\{ d,\, e,\, g \} \ \text{and}\ 
U_3=\{ b,\, e,\, g \}. \]
Each $U_i$ induces the separation ${\mathcal A}'_i$ of $S'$ 
as in Remark~\ref{Rem3-2}.
Then for any minimal separation ${\mathcal A}'$ of $S'$, 
${\mathcal A}' \preceq {\mathcal A}'_i$ for some $i=1,2,3$.
Here we note that the set $\{A_1,A_3,A_4\}$ is not a separation of $S'$.
Then the minimal separations of $S'$ are the 3-sets as 
\begin{align*}
&\{ A_1\cup A_3,\ A_4 \}, \\
&\{ A_1,\ A_3\cup A_4 \} \ \text{and} \\
&\{ A_1\cup A_4,\ A_3 \}.
\end{align*}
Hence, there are no type(I) subsets of $S'$ and 
\[ \widetilde{\mathcal A}_{S'}= \widetilde{\mathcal A}_{S'}^{\rm (II)}=\{S'\}. \]
\end{Example}

\begin{Example}\label{example10}
We consider Coxeter systems $(W,S)$ and $(W,S')$ 
defined by Figure~\ref{fig10}.
Here $S$ and $S'$ are the vertex sets of the corresponding figures and 
we define $o(a_1 a_2)=3$.

\begin{figure}[h]
\centering
{ 
\vspace{3mm}
\includegraphics[keepaspectratio, scale=0.95, bb=0 0 385 134]{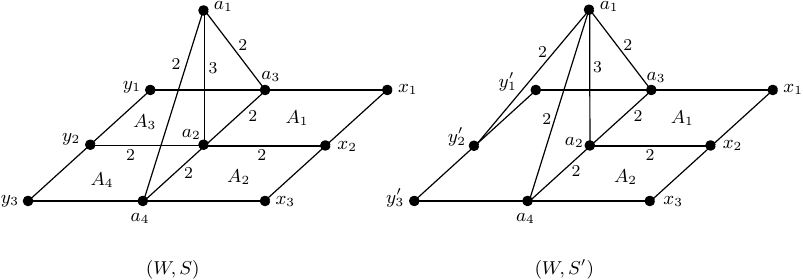}
}
\caption{Example~\ref{example10}}\label{fig10}
\end{figure}

Let $\sigma:=\{a_1,a_2\}$, $U:=\{a_1,a_2,a_3,a_4\}$, $X:=\{ x_1,x_2,x_3  \}$, $Y:=\{ y_1,y_2,y_3 \}$ 
and $Y':=\{ y'_1,y'_2,y'_3 \}$.
Here $U_\sigma=\sigma$, $U_\nu=\{a_3,a_4\}$ and $w_\sigma U=w_\sigma=U$.
Then $S'=X \cup U \cup w_\sigma Y w_\sigma$ is a twist of $S$ by $U$ and $w_\sigma$ 
(where $y'_i=w_\sigma y_i w_\sigma$ for $i=1,2,3,4$).

In $(W,S)$, 
let $A_i$ be the set of $4$-vertices around ``$A_i$'' in the figure for each $i=1,2,3,4$.
Then the minimal separations of $S$ are the 9-sets as 
\begin{enumerate}
\item[(1)] $\{U,\ A_1\cup A_2,\ A_3\cup A_4 \}$,
\item[(2)] $\{U,\ A_1,\ A_4,\ A_2\cup A_3 \}$,
\item[(3)] $\{U,\ A_1\cup A_4,\ A_2,\ A_3 \}$,
\item[(4)] $\{U\cup A_1,\  A_2,\ A_3\cup A_4 \}$,
\item[(5)] $\{U\cup A_2,\  A_1,\ A_3\cup A_4 \}$,
\item[(6)] $\{A_1\cup A_2,\  U\cup A_3,\ A_4 \}$,
\item[(7)] $\{A_1\cup A_2,\  U\cup A_4,\ A_3 \}$,
\item[(8)] $\{U\cup A_1\cup A_3,\  A_2,\ A_4 \}$ and
\item[(9)] $\{U\cup A_2\cup A_4,\  A_1,\ A_3 \}$.
\end{enumerate}

\smallskip

The Coxeter system $(W,S')$ is denoted by Figure~\ref{fig10-2}.
\begin{figure}[h]
\centering
{ 
\vspace{3mm}
\includegraphics[keepaspectratio, scale=0.95, bb=0 0 233 100]{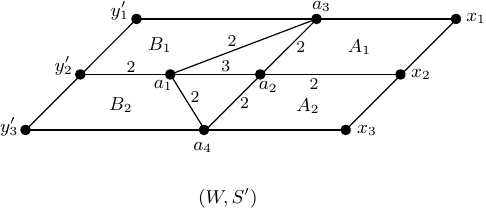}
}
\caption{The Coxeter system $(W,S')$}\label{fig10-2}
\end{figure}

In $(W,S')$, 
let $A_i$ and $B_i$ be the sets of $4$-vertices around ``$A_i$'' and ``$B_i$'' 
in the figure for each $i=1,2$ respectively.
Then the minimal separations of $S$ are the 9-sets as 
\begin{enumerate}
\item[(1)] $\{A_1\cup A_2,\ U,\ B_1\cup B_2 \}$,
\item[(2)] $\{A_1\cup A_2,\ U\cup B_1,\ B_2 \}$,
\item[(3)] $\{A_1\cup A_2,\ U\cup B_2,\ B_1 \}$,
\item[(4)] $\{A_1,\ A_2\cup U,\ B_1\cup B_2 \}$,
\item[(5)] $\{A_1,\ A_2\cup U\cup B_1,\ B_2 \}$,
\item[(6)] $\{A_1,\ A_2\cup U\cup B_2,\ B_1 \}$,
\item[(7)] $\{A_2,\ A_1\cup U,\ B_1\cup B_2 \}$,
\item[(8)] $\{A_2,\ A_1\cup U\cup B_1,\ B_2 \}$ and
\item[(9)] $\{A_2,\ A_1\cup U\cup B_2,\ B_1 \}$.
\end{enumerate}

Here $U$ is a spherical-product subset that separates $S$ and $S'$ both.
We consider the separations ${\mathcal A}$ and ${\mathcal A}'$ of $S$ and $S'$ 
induced by $U$ respectively.
Then $(W,S)$ and $(W,S')$ are compatible on the separations ${\mathcal A}$ and ${\mathcal A}'$.
Thus $(W,S)$ and $(W,S')$ are some-separation-compatible.

\smallskip

In $(W,S)$, there are no type(I) subsets of $S$ and 
\[ \widetilde{\mathcal A}_{S}= \widetilde{\mathcal A}_{S}^{\rm (II)}=\{S\}. \]

The Coxeter system $(W,S')$ is denoted by Figure~\ref{fig10-2}.
In $(W,S')$, there are no type(I) subsets of $S'$ and 
\[ \widetilde{\mathcal A}_{S'}= \widetilde{\mathcal A}_{S'}^{\rm (II)}=\{S'\}. \]

Then the type(II) subsets $S \in \widetilde{\mathcal A}_{S}^{\rm (II)}$ and 
$S' \in \widetilde{\mathcal A}_{S'}^{\rm (II)}$ are type(II)-compatible.
Hence $(W,S)$ and $(W,S')$ are type(I)-type(II)-compatible.
Here $(W,S)$ and $(W,S')$ are not compatible on their standard separations.
\end{Example}

\section{On separations of Coxeter generating sets and compatible Coxeter systems}\label{sec3}

Let $(W,S)$ be a Coxeter system.
A maximal spherical subset $\sigma$ of $S$ 
is not separated by any subset of $\sigma$.
Hence every maximal spherical subset $\sigma$ is a twist-rigid subset of $S$.
Also for a maximal spherical subset $\sigma$ of $S$, 
there exists a maximal twist-rigid subset $A$ of $S$ 
such that $\sigma \subset A$.

We consider that 
two Coxeter systems $(W,R)$ and $(W,S)$ 
are \textit{maximal-twist-rigid-subset-compatible}, 
if each maximal twist-rigid subset $A$ of $R$ 
is conjugate to some unique maximal twist-rigid subset $B$ of $S$ 
and each maximal twist-rigid subset $B$ of $S$ 
is conjugate to some unique maximal twist-rigid subset $A$ of $R$.

By the above argument, if two Coxeter systems are 
maximal-twist-rigid-subset-compatible, then 
they are maximal-spherical-subset-compatible.
Also by the argument in Section~\ref{sec1}, 
if two Coxeter systems are maximal-spherical-subset-compatible, 
then they are angle-compatible.

\begin{Question}
\begin{enumerate}
\item[(i)] Angle-compatible Coxeter systems $(W,R)$ and $(W,S)$ 
will be maximal-spherical-subset-compatible?
\item[(ii)] Maximal-spherical-subset-compatible Coxeter systems $(W,R)$ and $(W,S)$ 
will be maximal-twist-rigid-subset-compatible?
\end{enumerate}
\end{Question}

The following technical lemma and remark are 
used in the proof of the main theorem.

\begin{Lemma}\label{RemB1}
Let $(W,S)$ be a Coxeter system and 
let ${\mathcal A}$ be a separation of $S$.
Suppose that the following statements \textnormal{(1)--(4)} hold\textnormal{:}
\begin{enumerate}
\item[\textnormal{(1)}] $A_1,\ldots, A_n \in {\mathcal A}$.
\item[\textnormal{(2)}] $(A_1\cup \cdots \cup A_i)\cap A_{i+1}$ is maximal in 
\[ \{ (A_1\cup \cdots \cup A_i)\cap C : 
C \in {\mathcal A}-\{A_1,\ldots, A_i \} \} \]
for any $i=1,\ldots,n-1$.
\item[\textnormal{(3)}] $\overline{A}_0:=A_1\cup \cdots \cup A_n$.
\item[\textnormal{(4)}] $A \in {\mathcal A}-\{ A_1,\ldots, A_n \}$.
\end{enumerate}
Then there exists a subset $U_0$ of $\overline{A}_0$ such that 
$U_0$ is a separator of ${\mathcal A}$ and 
$U_0$ separates $\overline{A}_0$ and $A$; that is, 
for the connected components $X_1,\dots,X_t$ of $S-U_0$ 
(where $S-U_0=X_1 \cup \cdots \cup X_t$ is a disjoint union), 
$\overline{A}_0 \subset X_i \cup U_0$ and $A \subset X_j \cup U_0$ 
for some $i,j \in \{ 1,\ldots,t \}$ as $i\neq j$.
\end{Lemma}

\begin{proof}
We suppose that the statements (1)--(4) hold.

Let $A_{n+1} \in {\mathcal A}-\{ A_1,\ldots, A_n \}$ such that 
$U:=(A_1\cup \cdots \cup A_n)\cap A_{n+1}$ is maximal in 
$\{ (A_1\cup \cdots \cup A_n)\cap C : C \in {\mathcal A}-\{A_1,\ldots, A_n \} \}$.

If $U$ separates $\overline{A}_0$ and $A$, 
then we obtain $U_0:=U$ and this lemma is proved, 
since $U$ is a subset of $\overline{A}_0$ and 
$U$ is a separator of ${\mathcal A}$ by the definition of a separation.

We suppose that $U$ does not separate $\overline{A}_0$ and $A$.

Let $A_{n+1},\ldots,A_{p} \in {\mathcal A}-\{ A_1,\ldots, A_n \}$ such that 
$(A_1\cup \cdots \cup A_i)\cap A_{i+1}$ is maximal in 
\[ \{ (A_1\cup \cdots \cup A_i)\cap C : C \in {\mathcal A}-\{A_1,\ldots, A_i \} \} \]
for any $i=1,\ldots,n-1,n,\ldots,p-1$ and 
$(A_1\cup \cdots \cup A_{p-1})\cap A_{p}$ separates $\overline{A}_0$ and $A$.
Here we may suppose that $p$ is the minimum number 
as $(A_1\cup \cdots \cup A_{p-1})\cap A_{p}$ separates $\overline{A}_0$ and $A$;
that is, 
$(A_1\cup \cdots \cup A_{i-1})\cap A_{i}$ does not separate $\overline{A}_0$ and $A$ 
for any $i=n+1,\ldots,p-1$.

Let $U:=(A_1\cup \cdots \cup A_{p-1})\cap A_{p}$ 
that is a separator of ${\mathcal A}$ and separates $\overline{A}_0$ and $A$.
Here if $U \subset \overline{A}_0$ then we obtain $U_0:=U$ and this lemma is proved.

We suppose that $U \not\subset \overline{A}_0$.

We consider the connected components $Y_1,\dots,Y_r$ of $S-U$ 
where $S-U=Y_1 \cup \cdots \cup Y_r$ is a disjoint union.
We may suppose that $\overline{A}_0 \subset Y_1 \cup U$ and $A \subset Y_2 \cup U$.
By the definition of a separator, 
there exists $A_{i_0}$ ($1\le i_0 \le p-1$) such that $U \subset A_{i_0} \subset Y_1 \cup U$.
Since $U \not\subset \overline{A}_0$, we have that $n+1 \le i_0 \le p-1$.

\smallskip

(a) Since $U$ separates $\overline{A}_0$ and $A$, 
for any $a \in \overline{A}_0$ and $b \in A$,
any path $[a,b]$ form $a$ to $b$ in the nerve of $(W,S)$ intersects $U$.

\smallskip

(b) Let $V:=(A_1 \cup \cdots \cup A_{i_0-1})\cap A_{i_0}$ 
that is a separator of ${\mathcal A}$ and does not separate $\overline{A}_0$ and $A$.
We consider the connected components $Z_1,\dots,Z_q$ of $S-V$ 
where $S-V=Z_1 \cup \cdots \cup Z_q$ is a disjoint union.
We may suppose that $\overline{A}_0 \cup A \subset Z_1 \cup V$ and $A_{i_0} \subset Z_2 \cup V$.
Here $\overline{A}_0 \not\subset V$ and $A \not\subset V$.
Hence there exist $a_1 \in Z_1 \cap \overline{A}_0$ and $b_1 \in Z_1 \cap A$.
Since $Z_1$ is connected, 
there exists a path $[a_1,b_1]$ from $a_1$ to $b_1$ in the nerve of $(W_{Z_1},Z_1)$.
Here $Z_1 \cap A_{i_0} = \emptyset$, since $A_{i_0} \subset Z_2 \cup V$.
Thus the path $[a_1,b_1]$ does not intersect $A_{i_0}$.
In particular, the path $[a_1,b_1]$ does not intersect $U$, because $U \subset A_{i_0}$.
This contradicts (a).

\smallskip

Thus, $U \subset \overline{A}_0$ and we obtain $U_0:=U$. This lemma is proved.
\end{proof}

The following remark is obtained from Lemma~\ref{RemB1}.

\begin{Remark}\label{lemB}
Let $(W,S)$ be a Coxeter system with the untangle-conjugate-condition and 
let ${\mathcal A}$ be a separation of $S$.
We suppose that the following statements (1)--(8) hold\textnormal{:}
\begin{enumerate}
\item[\textnormal{(1)}] $A_1,\ldots, A_n \in {\mathcal A}$.
\item[\textnormal{(2)}] $(A_1\cup \cdots \cup A_i)\cap A_{i+1}$ is maximal in 
\[ \{ (A_1\cup \cdots \cup A_i)\cap C : 
C \in {\mathcal A}-\{A_1,\ldots, A_i \} \} \]
for any $i=1,\ldots,n-1$.
\item[\textnormal{(3)}] $\overline{A}_0:=A_1\cup \cdots \cup A_n$.
\item[\textnormal{(4)}] $A \in {\mathcal A}-\{ A_1,\ldots, A_n \}$.
\item[\textnormal{(5)}] $U$ is a spherical-product subset of $S$ 
such that $U_\sigma$ is non-empty and $U\subsetneq \overline{A}_0$.
\item[\textnormal{(6)}] $U'$ is a spherical-product subset of $S$ such that $U' \subsetneq A$.
\item[\textnormal{(7)}] $U$ and $U'$ are conjugate and $U\neq U'$.
\item[\textnormal{(8)}] $|U| \ge \max \{ \, |\overline{A}_0 \cap C| \, : 
C\in {\mathcal A},\ C \not\subset \overline{A}_0 \}$.
\end{enumerate}

Then by (7) and the untangle-conjugate-condition, 
there exist a sequence $U_1,\cdots, U_q$ of subsets of $S$ and 
a sequence $T_1,\cdots, T_{q-1}$ of spherical subsets of $S$ such that 
\[ \ds U=U_1 \mathop{\simeq}_{\,\,w_{T_1}} U_2 \mathop{\simeq}_{\,\,w_{T_2}} 
\cdots {\hspace*{-3mm}}\mathop{\simeq}_{\ \  w_{T_{q-1}}} U_q = U'. \]
Here $|U|=|U_1|=\cdots = |U_q|=|U'|$.

Since $U \neq U'$ by (7), $q \ge 2$.
We suppose that $U_i \not\subset \overline{A}_0$ for any $i \in \{ 2,\ldots,q \}$ 
and $T_i \not\subset \overline{A}_0$ for any $i \in \{ 1,\ldots,q-1 \}$.

By Lemma~\ref{RemB1}, 
there exists a subset $U_0$ of $\overline{A}_0$ such that 
$U_0$ is a separator of ${\mathcal A}$ and 
$U_0$ separates $\overline{A}_0$ and $A$.

By the above untangle-conjugate sequence, 
$F:=T_1 \cup \cdots \cup T_{q-1}$ is connected.
Here $(U_i)_{\sigma} \cup (U_{i+1})_{\sigma} \subset T_i$ for each $i=1,\ldots,q-1$.
Also $U_1=U\subsetneq \overline{A}_0$ by (5) and $U_q=U' \subsetneq A$ by (6).
Then $U \subset \overline{A}_0 \cap F$ and $U' \subset A \cap F$.
Also $U_{\nu}=U'_{\nu} \subset \overline{A}_0 \cap A$.

For any $a \in U_{\sigma}$, 
$F-(U-\{ a \})=(F-U)\cup \{ a \}$ is connected and there exists a path from $a$ to some point of $A$ in $F-(U-\{ a \})$, 
because each $T_{i}$ is spherical and $o(st) < \infty$ for any $s,t \in T_i$.

Hence $U \subset U_0$, since $U_0$ separates $\overline{A}_0$ and $A$.

We show that $U_0=\overline{A}_0 \cap C$ for some $C\in {\mathcal A}$ as $C \not\subset \overline{A}_0$.
Let $X_1,\ldots,X_t$ be the connected components of $S - U_0$ and 
let $\overline{X}_j:=\bigcup \{ C\in {\mathcal A} : C \subset X_j\cup U_0 \}$ 
for each $j=1,\ldots,t$.
Here $\overline{X}_j \cap \overline{X}_{j'} \subset U_0$ if $j\neq j'$.
Since $U_0$ is a separator of ${\mathcal A}$, 
there exist $C_1,C_2 \in {\mathcal A}$ such that 
$U_0=C_1 \cap C_2$, $C_1 \subset \overline{X}_{j_1}$ and $C_2 \subset \overline{X}_{j_2}$ 
for some $j_1,j_2 \in \{1,\ldots,t\}$ as $j_1 \neq j_2$ by (iv) in the definition of a separator.
Here $\overline{A}_0 \subset \overline{X}_{j_0}$ for some $j_0 \in \{1,\ldots,t\}$.
Then $j_0 \neq j_1$ or $j_0 \neq j_2$, since $j_1 \neq j_2$.
We suppose that $j_0 \neq j_1$.
Then $C_1 \not\subset \overline{A}_0$ and $U_0=\overline{A}_0 \cap C_1$, 
because $U_0 \subset \overline{A}_0$ and $U_0=C_1 \cap C_2 \subset C_1$ 
(here $\overline{X}_{j_0} \cap \overline{X}_{j_1} \subset U_0$).
Hence $U_0=\overline{A}_0 \cap C$ for some $C\in {\mathcal A}$ as $C \not\subset \overline{A}_0$.

Thus $|U| \ge |U_0|$ by (8).
Hence we have that $U= U_0$, because $U \subset U_0$.

Then $U$ separates $S$, and $V:=U \cup T_1=T_1 \cup U_{\nu}$ separates $S$.
Here $st=ts$ for any $s\in U_{\nu}$ and $t\in T_1$.

Let $X$ and $Y$ be subsets of $S$ such that 
$S-V =X \cup Y$ is a disjoint union, $\overline{A}_0 \subset X \cup V$, $A \subset Y \cup V$ 
and $o(xy)=\infty$ for any $x\in X$ and $y \in Y$.
Then we can obtain a twist 
\[ S'= X \cup V \cup w_{T_1} Y w_{T_1}. \]
Here for each $C \in {\mathcal A}$, 
$C \subset X \cup V$ or $C \subset Y \cup V$.
Hence the twist $S'$ is preserving ${\mathcal A}$ and fixing $\overline{A}_0$.
Let ${\mathcal A}'$ be the induced separation of $S'$ by ${\mathcal A}$ and the twist.
Then $w_{T_1} U_2 w_{T_1}=U_1=U$.
Let $A':=w_{T_1} A w_{T_1}$, 
let $U'_i:=w_{T_1} U_i w_{T_1}$ for $i=2,\ldots,q$ and 
let $T'_i:=w_{T_1} T_i w_{T_1}$ for $i=2,\ldots,q-1$.
Then in $S'$, we obtain that $U'_2 \subset \overline{A}_0$ and $U'_q \subset A'$.
Also 
\[ \ds U'_2 \mathop{\simeq}_{\,\,w_{T'_2}} U'_3 \mathop{\simeq}_{\,\,w_{T'_3}} 
\cdots {\hspace*{-3mm}}\mathop{\simeq}_{\ \  w_{T'_{q-1}}} U'_q \]
in $S'$.
Here we can iterate this argument.
\end{Remark}

\section{Some-separation-compatible Coxeter generating sets are conjugate up to finite twists}\label{sec4}

We prove that two some-separation-compatible Coxeter generating sets are conjugate up to finite twists.

\begin{Theorem}\label{maintheorem}
Let $(W,R)$ and $(W,S)$ be Coxeter systems with the untangle-condition.
If $(W,R)$ and $(W,S)$ are some-separation-compatible, then 
$R$ and $S$ are conjugate up to finite twists.
\end{Theorem}

\begin{proof}
We suppose that $(W,R)$ and $(W,S)$ are some-separation-compatible 
and suppose that $R$ and $S$ are connected.
Then there exist separations ${\mathcal A}_R$ and ${\mathcal A}_S$ 
of $R$ and $S$ respectively such that 
$(W,R)$ and $(W,S)$ are compatible on the separations ${\mathcal A}_R$ and ${\mathcal A}_S$; 
that is, 
each $A \in {\mathcal A}_{R}$ is conjugate to some unique $B \in {\mathcal A}_S$ and 
each $B \in {\mathcal A}_{S}$ is conjugate to some unique $A \in {\mathcal A}_R$.

We note that if $R'$ and $S'$ are Coxeter generating sets for $W$ 
that are obtained from $R$ and $S$ by some finite twists 
\textit{preserving} ${\mathcal A}_{R}$ and ${\mathcal A}_{S}$ respectively, 
then also 
each $A\in {\mathcal A}_{R'}$ is conjugate to some unique $B\in {\mathcal A}_{S'}$ and
each $B\in {\mathcal A}_{S'}$ is conjugate to some unique $A\in {\mathcal A}_{R'}$.

\medskip

Now we show that the following statement $({\rm P}_i)$ holds 
for each $i=1,\ldots,n$, 
where $n:=|{\mathcal A}_R|=|{\mathcal A}_S|$.
\begin{enumerate}
\item[$({\rm P}_i)$] 
There exist Coxeter generating sets $R_i$ and $S_i$ for $W$ such that 
$R_i$ and $S_i$ are obtained from $R$ and $S$ 
by some finite twists preserving ${\mathcal A}_{R}$ and ${\mathcal A}_S$ 
respectively and 
$\overline{A}_i$ and $\overline{B}_i$ are conjugate in $W$, 
where $\overline{A}_i$ is a union of 
some $i$-th elements $C_1,\ldots, C_i$ of ${\mathcal A}_{R_i}$ in $R_i$ 
as $(C_1\cup \cdots \cup C_k)\cap C_{k+1}$ is maximal in 
\[ \{ (C_1\cup \cdots \cup C_k)\cap C : 
C \in {\mathcal A}_{R_i} -\{C_1,\ldots, C_k \} \} \]
for any $k=1,\ldots,i-1$, 
and $\overline{B}_i$ is a union of 
some $i$-th elements $D_1,\ldots, D_i$ 
of ${\mathcal A}_{S_i}$ in $S_i$ 
as $(D_1\cup \cdots \cup D_k)\cap D_{k+1}$ is maximal in 
\[ \{ (D_1\cup \cdots \cup D_k)\cap D : 
D \in {\mathcal A}_{S_i} -\{D_1,\ldots, D_k \} \} \]
for any $k=1,\ldots,i-1$.
\end{enumerate}

We prove $({\rm P}_i)$ by induction on $i$.

Let $R_1:=R$ and $S_1:=S$.
For $A_1\in {\mathcal A}_{R_1}$, 
there exists a unique $B_1 \in {\mathcal A}_{S_1}$ 
such that $A_1$ and $B_1$ are conjugate.
Let $\overline{A}_1:=A_1$ and $\overline{B}_1:=B_1$.
Then $({\rm P}_1)$ holds.

Let $i_0\ge 1$.
We suppose that $({\rm P}_{i_0})$ holds; 
that is,
there exist Coxeter generating sets $R_{i_0}$ and $S_{i_0}$ 
for $W$ such that 
$R_{i_0}$ and $S_{i_0}$ are obtained from $R$ and $S$ 
by some finite twists preserving ${\mathcal A}_{R}$ and ${\mathcal A}_S$ 
respectively and 
$\overline{A}_{i_0}$ and $\overline{B}_{i_0}$ are conjugate, 
where $\overline{A}_{i_0}$ is a union of 
some $i_0$-th elements $C_1,\ldots, C_{i_0}$ 
of ${\mathcal A}_{R_{i_0}}$ in $R_{i_0}$ 
as $(C_1\cup \cdots \cup C_k)\cap C_{k+1}$ is maximal in 
\[ \{ (C_1\cup \cdots \cup C_k)\cap C : 
C \in {\mathcal A}_{R_{i_0}} -\{C_1,\ldots, C_k \} \} \]
for any $k=1,\ldots,i_0-1$, 
and $\overline{B}_{i_0}$ is a union of 
some $i_0$-th elements $D_1,\ldots, D_{i_0}$ 
of ${\mathcal A}_{S_{i_0}}$ in $S_{i_0}$ 
as $(D_1\cup \cdots \cup D_k)\cap D_{k+1}$ is maximal in 
\[ \{ (D_1\cup \cdots \cup D_k)\cap D : 
D \in {\mathcal A}_{S_{i_0}} -\{D_1,\ldots, D_k \} \} \]
for any $k=1,\ldots,i_0-1$.

Then we show that $({\rm P}_{i_0 +1})$ holds.

Let 
\begin{align*}
&d_1:= \max \{\; |\overline{A}_{i_0} \cap C| \, : 
C\in {\mathcal A}_{R_{i_0}}, \ 
C \not\subset \overline{A}_{i_0} \} \text{ and}\\
&d_2:= \max \{\; |\overline{B}_{i_0} \cap D| \, : 
D\in {\mathcal A}_{S_{i_0}}, \ 
D \not\subset \overline{B}_{i_0} \}.
\end{align*}

Since $R$ and $S$ are connected, 
$R_{i_0}$ and $S_{i_0}$ are connected.
Hence $d_1>0$ and $d_2>0$.

Now we suppose that $d_1 \le d_2$.

Then there exists $B\in {\mathcal A}_{S_{i_0}}$ 
such that $B \not\subset \overline{B}_{i_0}$ and 
$d_2= |\overline{B}_{i_0} \cap B|$.
Let $S_{i_0+1}:=S_{i_0}$ and 
let $\overline{B}_{i_0+1}:= \overline{B}_{i_0} \cup B$.
Here $S_{i_0+1}$ is preserving ${\mathcal A}_{S_{i_0}}$ and fixing $\overline{B}_{i_0}$.
Also $S_{i_0+1}$ is obtained from $S$ by some finite twists preserving ${\mathcal A}_S$.
We also note that $\overline{B}_{i_0+1}$ is a union 
of some $(i_0+1)$-th elements of ${\mathcal A}_{S_{i_0+1}}$ in $S_{i_0+1}$.

Let $A$ be the element of ${\mathcal A}_{R_{i_0}}$ 
such that $A$ and $B$ are conjugate in $W$.
Let $U:=\overline{B}_{i_0} \cap B$.
Then since $|U|=d_2$, 
$U$ is a spherical-product subset of $S_{i_0+1}$ 
and $U$ separates $S_{i_0+1}$ 
by the definition of a separation, 
because ${\mathcal A}_{S_{i_0 +1}}$ is a separation of $S_{i_0 +1}$.

Since $\overline{A}_{i_0}$ and $\overline{B}_{i_0}$ are conjugate in $W$, 
$x \overline{B}_{i_0} x^{-1}= \overline{A}_{i_0}$ for some $x\in W$.
Then $U \subset \overline{B}_{i_0}$ and 
$U':=x U x^{-1} \subset \overline{A}_{i_0}$.
Also since $A$ and $B$ are conjugate in $W$, 
$y B y^{-1}= A$ for some $y\in W$.
Then $U \subset B$ and 
$U'':=y U y^{-1} \subset A$.
Here $U'$ and $U''$ are 
conjugate spherical-product subsets of $R_{i_0}$ in $W$.
Hence they are untangle by the untangle-condition.
Here if $U'_\sigma$ is empty then 
$U'=U'_\nu =U''_\nu =U''$.

\medskip

(a) We consider the case that $U'=U''$.

Then $\overline{A}_{i_0} \cap A=U'=U''$ and $d_1=d_2$, 
since $|U'|=|U|=d_2$ and $d_1 \le d_2$.
Hence $|U'|=d_1$, 
$U'=\overline{A}_{i_0}\cap A$ is maximal in 
\[ \{ \overline{A}_{i_0}\cap C : 
C \in {\mathcal A}_{R_{i_0}}, \ C\not\subset \overline{A}_{i_0} \} \]
and $U'=U''$ separates $R_{i_0}$.

Here $x U x^{-1}=U'$ and $y U y^{-1}=U''=U'$.
Hence 
\[ x y^{-1} U' y x^{-1} = x U x^{-1} =U'. \]

We define the bijective map $f_{x y^{-1}}:U' \to U'$ 
by $f_{x y^{-1}}(a)=(x y^{-1}) a (y x^{-1})$ for any $a \in U'$.

By the untangle-conjugate-condition, 
$f_{x y^{-1}}=f_1$ for $1\in W$, or, 
there exist a sequence $U_1,\cdots, U_q$ of subsets of $R_{i_0}$ and 
a sequence $T_1,\cdots, T_{q-1}$ of spherical subsets of $R_{i_0}$ such that 
\[ \ds U'=U_1 \mathop{\simeq}_{\,\,w_{T_1}} U_2 \mathop{\simeq}_{\,\,w_{T_2}} 
\cdots {\hspace*{-3mm}}\mathop{\simeq}_{\ \  w_{T_{q-1}}} U_q = U', \]
$(U_i)_{\sigma} \cup (U_{i+1})_{\sigma} \subset T_i$ for any $i=1,\ldots,q-1$ 
and $f_{x y^{-1}} = f_{w_0}$ for $w_0:=w_{T_{q-1}}\cdots w_{T_2}w_{T_1}$.

\medskip

(a-1) Suppose that $f_{x y^{-1}}=f_1$.
Then $\overline{A}_{i_0+1}=\overline{A}_{i_0} \cup A$ in $R_{i_0}$ 
is conjugate to $\overline{B}_{i_0+1}=\overline{B}_{i_0} \cup B$ in $S_{i_0}$ 
where $\overline{A}_{i_0+1}$ and $\overline{B}_{i_0+1}$ are 
unions of some $(i_0+1)$-th elements of ${\mathcal A}_{R_{i_0}}$ 
and ${\mathcal A}_{S_{i_0}}$ respectively.
Let $R_{i_0+1}:=R_{i_0}$ and $S_{i_0+1}:=S_{i_0}$.
Here $\overline{A}_{i_0} \cap A$ is maximal in 
\[ \{ \overline{A}_{i_0}\cap C : 
C \in {\mathcal A}_{R_{i_0+1}}, \ C\not\subset \overline{A}_{i_0} \} \]
and $\overline{B}_{i_0} \cap B$ is maximal in 
\[ \{ \overline{B}_{i_0}\cap D : 
D \in {\mathcal A}_{S_{i_0+1}}, \ D\not\subset \overline{B}_{i_0} \}, \]
because $|\overline{A}_{i_0} \cap A|=|\overline{B}_{i_0} \cap B|=d_1=d_2$.

\medskip

(a-2) Suppose that there exist a sequence $U_1,\cdots, U_q$ of subsets of $R_{i_0}$ and 
a sequence $T_1,\cdots, T_{q-1}$ of spherical subsets of $R_{i_0}$ such that 
\[ \ds U'=U_1 \mathop{\simeq}_{\,\,w_{T_1}} U_2 \mathop{\simeq}_{\,\,w_{T_2}} 
\cdots {\hspace*{-3mm}}\mathop{\simeq}_{\ \  w_{T_{q-1}}} U_q = U', \]
$(U_i)_{\sigma} \cup (U_{i+1})_{\sigma} \subset T_i$ for any $i=1,\ldots,q-1$ 
and $f_{x y^{-1}} = f_{w_0}$ for $w_0:=w_{T_{q-1}}\cdots w_{T_2}w_{T_1}$.

Then we will attach $A$ to $\overline{A}_{i_0}$ by gluing $U'$ and $U''=U'$ 
by some finite twists of $R_{i_0}$ 
(preserving ${\mathcal A}_{R_{i_0}}$ and fixing $\overline{A}_{i_0}$) 
induced by the above untangle-conjugate sequence from $U'$ to $U''=U'$.

Since $U'=\overline{A}_{i_0}\cap A$ separates $R_{i_0}$, 
$V:=U' \cup T_{1}=T_{1} \cup U'_{\nu}$ separates $R_{i_0}$.
For some two subsets $X$ and $Y$ of $R_{i_0}$, 
\begin{enumerate}
\item[(1)] $R_{i_0}-V =X \cup Y$ that is a disjoint union, 
\item[(2)] $o(xy)=\infty$ for any $x\in X$ and $y\in Y$, 
\item[(3)] $\overline{A}_{i_0} \cup (U_1\cup\cdots\cup U_q) \subset X\cup V$ and 
\item[(4)] $A \subset Y\cup V$.
\end{enumerate}
Here for any $C\in {\mathcal A}_{R_{i_0}}$, 
$C \subset X\cup V$ or $C \subset Y\cup V$, 
since $U'=\overline{A}_{i_0} \cap A$ is a separator of ${\mathcal A}_{R_{i_0}}$ and $V=U' \cup T_{1}$.

Then $R'_{i_0}:= X \cup V \cup (w_{T_{1}} Y w_{T_{1}})$ is a twist of $R_{i_0}$ 
preserving ${\mathcal A}_{R_{i_0}}$ and fixing $\overline{A}_{i_0}$.
Let ${\mathcal A}_{R'_{i_0}}$ be the separation of $R'_{i_0}$ induced by 
${\mathcal A}_{R_{i_0}}$ and the twist.

Let $A':=w_{T_{1}} A w_{T_{1}} \in {\mathcal A}_{R'_{i_0}}$ that is conjugate to $A$ and $B$.
Then 
\[ A'=w_{T_{1}} A w_{T_{1}} \supset w_{T_{1}} U' w_{T_{1}} = w_{T_{1}} U_1 w_{T_{1}} = U_{2}, \]
$U'=U''=U_q \subset \overline{A}_{i_0}$ and 
\[ \ds U_2 \mathop{\simeq}_{\,\,w_{T_2}} U_3 \mathop{\simeq}_{\,\,w_{T_3}} 
\cdots {\hspace*{-3mm}}\mathop{\simeq}_{\ \  w_{T_{q-1}}} U_{q}=U'. \]
Since $U'=U_1$ separates $R_{i_0}$ and separates $\overline{A}_{i_0}$ and $A$, 
we have that $U_2$ separates $R'_{i_0}$ and $U_2$ separates $\overline{A}_{i_0}$ and $A'$ in $R'_{i_0}$.
Hence $T_2 \cup U'_{\nu}$ separates $R'_{i_0}$ and $T_2 \cup U'_{\nu}$ separates $\overline{A}_{i_0}$ and $A'$ in $R'_{i_0}$.

We iterate this argument for $T_2,\ldots,T_{q-1}$.
Then we obtain a Coxeter generating set $R''_{i_0}$ from $R_{i_0}$ by some finite twists 
preserving ${\mathcal A}_{R_{i_0}}$ and fixing $\overline{A}_{i_0}$ 
such that 
for $w_0=w_{T_{q-1}}\cdots w_{T_1}$ and $A'':= w_0 A w_0 \in {\mathcal A}_{R''_{i_0}}$ that is conjugate to $A$ and $B$, 
\[ A''=w_0 A w_0 \supset w_0 U' w_0 = w_{T_{q-1}} \cdots w_{T_1} U_1 w_{T_1}\cdots w_{T_{q-1}} 
= U_q =U', \]
where ${\mathcal A}_{R''_{i_0}}$ is the separation of $R''_{i_0}$ induced 
by ${\mathcal A}_{R_{i_0}}$ and the twists.

Let $R_{i_0+1}:=R''_{i_0}$ and $S_{i_0+1}:=S_{i_0}$.
Here $\overline{A}_{i_0} \cap A''$ is maximal in 
\[ \{ \overline{A}_{i_0}\cap C : 
C \in {\mathcal A}_{R_{i_0+1}}, \ C\not\subset \overline{A}_{i_0} \} \]
and $\overline{B}_{i_0} \cap B$ is maximal in 
\[ \{ \overline{B}_{i_0}\cap D : 
D \in {\mathcal A}_{S_{i_0+1}}, \ D\not\subset \overline{B}_{i_0} \}, \]
because $|\overline{A}_{i_0} \cap A''|=|\overline{B}_{i_0} \cap B|=d_1=d_2$.
Here $\overline{A}_{i_0} \cup A''$ in $R_{i_0+1}$ is conjugate to 
$\overline{B}_{i_0} \cup B$ in $S_{i_0+1}$ 
where $A''= w_0 A w_0$ is the element of ${\mathcal A}_{R_{i_0+1}}$ 
that is conjugate to $A$ and $B$.
Thus we obtain that $({\rm P}_{i_0+1})$ holds.

\medskip

(b) We consider the case that $U' \neq U''$.

Here $U'_\sigma$ and $U''_\sigma$ are non-empty.
Then there exists a sequence of spherical-product subsets 
\[ \ds U' =U_1 \mathop{\simeq}_{\,\,w_{T_1}} U_2 \mathop{\simeq}_{\,\,w_{T_2}} \cdots 
{\hspace*{-3mm}}\mathop{\simeq}_{\ \ w_{T_{q-1}}} U_q=U'' \]
in $R_{i_0}$.
Here $T_1,\cdots, T_{q-1}$ are spherical subsets of $R_{i_0}$ such that 
$(U_j)_{\sigma} \cup (U_{j+1})_{\sigma} \subset T_j$ for $j=1,\ldots,q-1$, 
$(U_j)_\nu=U'_\nu$ for $j=1,\ldots,q$, 
and $st=ts$ for any $s\in U'_\nu$ and $t \in T_1\cup \cdots \cup T_{q-1}$.
Also $U'=U_1 \subset \overline{A}_{i_0}$ and $U''=U_q \subset A$.

Let $\sigma_j:=(U_j)_\sigma$ for each $j=1,\ldots,q$.
Then 
\[ \ds U'_\sigma =\sigma_1 \mathop{\simeq}_{\,\,w_{T_1}} \sigma_2 \mathop{\simeq}_{\,\,w_{T_2}} \cdots 
{\hspace*{-3mm}}\mathop{\simeq}_{\ \  w_{T_{q-1}}} \sigma_q = U''_\sigma. \]

We will attach $A$ to $\overline{A}_{i_0}$ by gluing $U'$ and $U''$ 
by some finite twists of $R_{i_0}$ 
(preserving ${\mathcal A}_{R_{i_0}}$ and fixing $\overline{A}_{i_0}$) 
induced by the untangle-conjugate sequence from $U'$ to $U''$.

\smallskip

We first suppose that $U_q \not\subset \overline{A}_{i_0}$ 
(that is, $\sigma_q \not\subset \overline{A}_{i_0}$).
Let $j_0 \in\{ 1,\ldots,q-1 \}$ be the number 
as $\sigma_{j_0} \subset \overline{A}_{i_0}$ 
and $\sigma_j \not\subset \overline{A}_{i_0}$ for any $j=j_0+1,\ldots,q$.
For $T:=T_{j_0}$, 
$w_T \sigma_{j_0+1} w_T = \sigma_{j_0}$ and 
$T \not\subset \overline{A}_{i_0}$.

Let $V:=T \cup U'_\nu=U_{j_0} \cup T$.
Here $st=ts$ for any $s\in U'_\nu$ and $t\in T$.
Then $T$ is a spherical subset 
and $V$ is a spherical-product subset of $R_{i_0}$.
Also $\overline{A}_{i_0} \cap V = U_{j_0}$ and 
$|U_{j_0}|=|U'|=|U|=d_2 \ge d_1$.

Then by Remark~\ref{lemB}, $U_{j_0}$ separates $R_{i_0}$ and 
$U_{j_0}$ separates $\overline{A}_{i_0}$ and $A$ in $R_{i_0}$.

Thus $V=U_{j_0}\cup T$ separates $R_{i_0}$ and 
$V$ separates $\overline{A}_{i_0}$ and $A$.
For some two subsets $X$ and $Y$ of $R_{i_0}$, 
\begin{enumerate}
\item[(1)] $R_{i_0}-V =X \cup Y$ that is a disjoint union, 
\item[(2)] $o(xy)=\infty$ for any $x\in X$ and $y\in Y$, 
\item[(3)] $\overline{A}_{i_0} \subset X\cup V$ and 
\item[(4)] $A \subset Y\cup V$.
\end{enumerate}
Here for any $C\in {\mathcal A}_{R_{i_0}}$, 
$C \subset X\cup V$ or $C \subset Y\cup V$, 
since $U_{j_0}=\overline{A}_{i_0} \cap A''$ and 
$V=U_{j_0} \cup T$.
Then $R'_{i_0}:= X \cup V \cup (w_T Y w_T)$ 
is a twist of $R_{i_0}$ 
preserving ${\mathcal A}_{R_{i_0}}$ and fixing $\overline{A}_{i_0}$.
Let ${\mathcal A}_{R'_{i_0}}$ be the separation of $R'_{i_0}$ induced by 
${\mathcal A}_{R_{i_0}}$ and the twist.

Then $w_T \sigma_{j_0+1} w_T=\sigma_{j_0}$ and 
$w_T U_{j_0+1} w_T=U_{j_0}$.
Hence $U''$ moves one step toward $U'$ by this twist.
Let $\sigma'_k:= w_T \sigma_k w_T$ and 
$U'_k:= w_T U_k w_T$ for each $k=j_0+1,\ldots,q$ 
and let $T'_k:= w_T T_k w_T$ for $k=j_0+1,\ldots,q-1$ 
(that are the corresponding subsets of $R'_{i_0}$ 
to $\sigma_k$, $U_k$ and $T_k$ in $R_{i_0}$ respectively).
Then 
\[ \ds U'_{\sigma} = \sigma_1 \mathop{\simeq}_{\,\,w_{T_1}} \cdots 
{\hspace*{-3mm}}\mathop{\simeq}_{\ \  w_{T_{j_0-1}}} \sigma_{j_0} =\sigma'_{j_0+1} 
{\hspace*{-3mm}}\mathop{\simeq}_{\ \  w_{T'_{j_0+1}}} \sigma'_{j_0+2} {\hspace*{-3mm}}\mathop{\simeq}_{\ \  w_{T'_{j_0+2}}} \cdots 
{\hspace*{-3mm}}\mathop{\simeq}_{\ \  w_{T'_{q-1}}} \sigma'_{q}=(\overline{U}'')_\sigma \]
and 
\[ \ds U'=U_1 \mathop{\simeq}_{\,\,w_{T_1}} \cdots {\hspace*{-3mm}}\mathop{\simeq}_{\ \  w_{T_{j_0-1}}} U_{j_0} = U'_{j_0+1} 
{\hspace*{-3mm}}\mathop{\simeq}_{\ \  w_{T'_{j_0+1}}} U'_{j_0+2} {\hspace*{-3mm}}\mathop{\simeq}_{\ \  w_{T'_{j_0+2}}} \cdots 
{\hspace*{-3mm}}\mathop{\simeq}_{\ \  w_{T'_{q-1}}} U'_{q} = \overline{U}'' \]
in $R'_{i_0}$ where $\overline{U}'':=w_T U'' w_T$.
Here $A':=w_T A w_T \in {\mathcal A}_{R'_{i_0}}$ is conjugate to $A$ and $B$, 
and $\overline{U}'' \subset A'$ holds.

We iterate this argument.
Then we obtain a Coxeter generating set $R''_{i_0}$ from $R_{i_0}$ 
by some finite twists 
preserving ${\mathcal A}_{R_{i_0}}$ and fixing $\overline{A}_{i_0}$ 
such that 
the above sequence of spherical-product subsets 
transforms to 
\[ \ds U'=U_1 \mathop{\simeq}_{\,\,w_{T_1}} U_2 
\mathop{\simeq}_{\,\,w_{T_2}} \cdots {\hspace*{-3mm}}\mathop{\simeq}_{\ \  w_{T_{j_0-1}}} U_{j_0}=U''' \]
in $R''_{i_0}$ where 
${\mathcal A}_{R''_{i_0}}$ is the separation of $R''_{i_0}$ induced by 
${\mathcal A}_{R_{i_0}}$ and the twist, and 
$U''' \subset A''$ for $A''\in {\mathcal A}_{R''_{i_0}}$ that is conjugate to $A$ and $B$.
Here $A'' \not\subset \overline{A}_{i_0}$ and 
$U_{j_0} \subset \overline{A}_{i_0}$ 
by the assumption and the definition of the number $j_0$.

\medskip

Then $U_{j_0}=\overline{A}_{i_0} \cap A''$ and $|U_{j_0}|=d_1=d_2$.
Hence $U_{j_0}$ separates $R''_{i_0}$.
Thus if $j_0 \ge 2$ then for $T_0:=T_{j_0-1}$, 
$V_0:=T_0 \cup U'_\nu = U_{j_0} \cup T_0$ separates $R''_{i_0}$ and 
we have a twist $R'''_{i_0}$ of $R''_{i_0}$ by $V_0$ and $w_{T_0}$ 
preserving ${\mathcal A}_{R''_{i_0}}$ and fixing $\overline{A}_{i_0}$.
Then $U'''$ moves one step toward $U'$ by this twist.
Here $A''':= w_{T_0} A'' w_{T_0}$ is conjugate to $A$ and $B$, and 
$U_{j_0-1}=w_{T_0} U_{j_0} w_{T_0} \subset A'''$ in $R'''_{i_0}$.

Since $U_{j_0}$ separates $R''_{i_0}$ and separates $\overline{A}_{i_0}$ and $A''$ in $R''_{i_0}$, 
we have that $U_{j_0-1}$ separates $R'''_{i_0}$ and $U_{j_0-1}$ separates $\overline{A}_{i_0}$ and $A'''$ in $R'''_{i_0}$.
Hence $T_{j_0-1} \cup U'_{\nu}$ separates $R'''_{i_0}$ and 
$T_{j_0-1} \cup U'_{\nu}$ separates $\overline{A}_{i_0}$ and $A'''$ in $R'''_{i_0}$.

We can iterate this argument for $T_{j_0-1},\ldots,T_1$.

\medskip

By iterating this argument and by (a), 
we obtain a Coxeter generating set $R_{i_0+1}$ from $R_{i_0}$ by some finite twists 
preserving ${\mathcal A}_{R_{i_0}}$ and fixing $\overline{A}_{i_0}$ 
such that 
$\overline{A}_{i_0} \cup A'$ in $R_{i_0+1}$ is conjugate to 
$\overline{B}_{i_0} \cup B$ in $S_{i_0+1}$ 
where $A'$ is the element of ${\mathcal A}_{R_{i_0+1}}$ 
that is conjugate to $A$ and $B$.
Hence, we obtain that $({\rm P}_{i_0+1})$ holds.

\medskip

Thus $({\rm P}_{i})$ holds for any $i=1,\ldots,n$.

Then $({\rm P}_{n})$ implies that 
there exist Coxeter generating sets $R_n$ and $S_n$ for $W$ such that 
$R_n$ and $S_n$ are obtained from $R$ and $S$ 
by some finite twists respectively and 
$R_n$ and $S_n$ are conjugate, 
because 
the unions of the $n$-th elements of ${\mathcal A}_{R_n}$ in $R_n$ 
and ${\mathcal A}_{S_n}$ in $S_n$ 
are just $R_n$ and $S_n$ respectively.

Therefore, $R$ and $S$ are conjugate up to finite twists.

\medskip

In the case that $R$ and $S$ are connected, 
we showed that 
if $(W,R)$ and $(W,S)$ are some-separation-compatible, 
then $R$ and $S$ are conjugate up to finite twists.

\smallskip

We suppose that $R$ and $S$ are not connected.

Let $R_1,\ldots,R_n$ be the connected components of $R$ and 
let $S_1,\ldots,S_n$ be the connected components of $S$; 
that is, 
\[ W=W_{S_1}*\cdots *W_{S_n}=W_{R_1}*\cdots *W_{R_n}, \]
where the numbers of the connected components of $R$ and $S$ are equal.

By the same argument as above (in the case that $R$ and $S$ are connected), 
we can obtain Coxeter generating sets $R'$ and $S'$ for $W$ such that 
$R'$ and $S'$ are obtained from $R$ and $S$ 
by some finite twists respectively and 
$R'_i$ and $S'_i$ are conjugate for any $i=1,\ldots,n$, 
where $R'_1,\ldots,R'_n$ are the connected components of $R'$ and 
$S'_1,\ldots,S'_n$ are the connected components of $S'$.

Then by the untangle-condition, 
there exists a Coxeter generating set $R''$ for $W$ 
such that $R''$ is obtained from $R'$ by some finite twists and $S'$ and $R''$ are conjugate.

Therefore, $R$ and $S$ are conjugate up to finite twists.
\end{proof}

We obtain the following corollary from Theorem~\ref{maintheorem}.

\begin{Corollary}\label{cor}
For Coxeter systems $(W,R)$ and $(W,S)$ with the untangle-condition, 
the following two statements are equivalent\textnormal{:}
\begin{enumerate}
\item[\textnormal{(i)}] $R$ and $S$ are conjugate up to finite twists.
\item[\textnormal{(ii)}] $(W,R_0)$ and $(W,S)$ are some-separation-compatible 
for some Coxeter generating set $R_0$ obtained from $R$ by finite twists.
\end{enumerate}
\end{Corollary}

\section{$\widetilde{\mathcal A}_S$ is a separation of $S$}\label{sec5}

In this section, we prove the following proposition.

\begin{Proposition}
Let $(W,S)$ be a Coxeter system 
and let $\widetilde{\mathcal A}_S$ be the set of subsets of $S$ 
in Definition~\ref{def:standard-separation}.
Then $\widetilde{\mathcal A}_S$ is a separation of $S$.
\end{Proposition}

Let $(W,S)$ be a Coxeter system 
and let ${\mathcal A}_0$ be the set of maximal twist-rigid subsets of $S$.
We suppose that $S$ is connected.
For $\widetilde{\mathcal A}_S$, we show the statements (1)--(4) 
of the definition of a separation (in Section~\ref{sec1}) hold.

\smallskip

(1) and (3): 
For each $A_0 \in {\mathcal A}_0$, there exists a unique $B \in \widetilde{\mathcal A}_S$ 
such that $A_0 \subset B$.
(Here we can denote $B=\bigcup [A_0]$ by Definition~\ref{def:standard-separation}.)
Hence $\bigcup \widetilde{\mathcal A}_S = \bigcup {\mathcal A}_0=S$.

\smallskip

(2) Let $B\in \widetilde{\mathcal A}_S$.
Then $B=\bigcup [A_0]$ for some $A_0\in {\mathcal A}_0$ 
by Definition~\ref{def:standard-separation}.
Hence $B$ is a union of some maximal twist-rigid subsets of $S$.
Also $B=\bigcup [A_0]$ is connected 
by the definition of the equivalence relation ``$\sim$'' on ${\mathcal A}_0$.

\smallskip

To prove (4) in the definition of a separation for $\widetilde{\mathcal A}_S$, 
we show some lemmas.

\smallskip

Let ${\mathcal A}$ be a minimal separation of $S$ (and fix ${\mathcal A}$).

\begin{Lemma}\label{lemmaA1}
Let $B_0\in \widetilde{\mathcal A}_S^{\rm (II)}$.
Let $A_1,\ldots,A_l$ be the distinct elements of ${\mathcal A}$ 
such that $A_i \subset B_0$ for any $i=1,\ldots,l$.
Then $\textstyle B_0=\bigcup_{i=1}^l A_i$.
\end{Lemma}

\begin{proof}
We can denote $B_0=\bigcup [A']$ for some $A' \in {\mathcal A}_0$.
Let $[A']=\{ A'_1,\ldots,A'_t \}$ (where $A'_j\in {\mathcal A}_0$ as $A'_j \sim A'$).
For each $j=1,\ldots,t$, 
$A'_j \subset A''_j$ for some unique $A''_j \in {\mathcal A}$, 
since ${\mathcal A}$ is a separation of $S$.
Here $A''_1,\ldots, A''_t$ need not be different all together.

We show that $A''_j\subset B_0$ for any $j \in \{1,\ldots,t\}$.

Let $j \in \{1,\ldots,t\}$.
Here $A'_j \subset B_0=\bigcup [A']$ and $A'_j \in {\mathcal A}_0$.
Hence $A'_j \sim A'$.
Let $A'' \in {\mathcal A}_0$ as $A'' \subset A''_j$.
Then $A'' \sim A'_j$ by the definition of the equivalence relation ``$\sim$'' on ${\mathcal A}_0$, 
since $A'',A'_j \in {\mathcal A}_0$ and $A''\cup A'_j \subset A''_j\in {\mathcal A}$.
Thus $A'' \sim A'_j \sim A'$.
We obtain that $A'' \subset \bigcup [A']=B_0$.
Hence 
\[ A''_j = \bigcup\{A'' \in {\mathcal A}_0 : A'' \subset A''_j \}\subset B_0. \]
Thus $A''_j \subset B_0$ for any $j \in \{1,\ldots,t\}$.

Then $\{A''_1,\ldots,A''_t\} \subset \{A_1,\ldots,A_l\}$ (in fact, the equality holds) 
and 
\[ \textstyle B_0=\bigcup [A']=\bigcup_{j=1}^t A'_j \subset \bigcup_{j=1}^t A''_j \subset \bigcup_{i=1}^l A_i. \]
Also obviously $\bigcup_{i=1}^l A_i \subset B_0$ holds.
Thus $B_0 = \bigcup_{i=1}^l A_i$.
\end{proof}

Now we say that a sequence $A_1,\ldots, A_l \in {\mathcal A}$ satisfies 
the condition $(*)$ in ${\mathcal A}$, if
\begin{enumerate}
\item[$(*)$\ \ ] 
$(A_{1}\cup \cdots \cup A_{i}) \cap A_{i+1}$ is 
maximal in 
\[ \{(A_{1}\cup \cdots \cup A_{i}) \cap A : 
A \in {\mathcal A}-\{A_{1},\ldots,A_{i}\} \} \]
for any $i=1,\ldots,l-1$.
\end{enumerate}

\medskip

We next show the following.

\begin{Lemma}\label{lemmaA2}
Let $B_0\in \widetilde{\mathcal A}_S^{\rm (II)}$.
Let $A_1,\ldots,A_l$ be the distinct elements of ${\mathcal A}$ 
such that $A_i \subset B_0$ for any $i=1,\ldots,l$.
Then there exists a bijective map $f:\{1,\ldots,l\} \to \{1,\ldots,l\}$ such that 
the sequence $A_{f(1)},\ldots,A_{f(l)}$ 
satisfies the condition $(*)$ in ${\mathcal A}$.
\end{Lemma}

\begin{proof}
Let $f(1)\in \{ 1,\ldots,l \}$ that is arbitrary.
Let $f(2)\in \{ 1,\ldots,l \}-\{f(1)\}$ such that 
$A_{f(1)} \cap A_{f(2)}$ is maximal in 
$\{A_{f(1)} \cap A : A\in \{ A_1,\ldots,A_l \}-\{A_{f(1)}\} \}$.
Then we show that $A_{f(1)} \cap A_{f(2)}$ is maximal in 
$\{A_{f(1)} \cap A : A\in {\mathcal A}-\{A_{f(1)}\} \}$.
Indeed if this does not hold, then 
there exists $A'_2\in {\mathcal A}-\{A_{f(1)}\}$ such that 
$A_{f(1)} \cap A'_2$ is maximal in 
$\{A_{f(1)} \cap A : A\in {\mathcal A}-\{A_{f(1)}\} \}$ 
and $A_{f(1)} \cap A_{f(2)} \subsetneq A_{f(1)} \cap A'_2$.
Here $A'_2 \not\in \{A_1,\ldots,A_l \}$ and $A'_2 \not\subset B_0$.
Hence some $\overline{U}_1 \in \overline{\mathcal U}$ separates $B_0$ and $A'_2$, 
where $\overline{\mathcal U}$ is the separators set as in Definition~\ref{def:standard-separation}.
Here ``$\overline{U}_1$ separates $B_0$ and $A'_2$'' means that 
$B_0 \subset \overline{X}_{j_1}$ and $A'_2 \subset \overline{X}_{j_2}$ 
for some $j_1,j_2\in\{1,\ldots,t\}$ as $j_1\neq j_2$, 
where $X_1,\ldots,X_t$ are the connected components of $S-U_1$ and 
$S-U_1=X_1 \cup \cdots \cup X_t$.

Then $\overline{U}_1$ separates $A_{f(1)}$ and $A'_2$.
Thus $\overline{U}_1$ separates $A_{f(1)}$ and $A_{f(2)}$, 
because $A_{f(1)} \cap A_{f(2)} \subsetneq A_{f(1)} \cap A'_2 \subset \overline{U}_1$.
This contradicts that 
$A_{f(1)} \cup A_{f(2)} \subset B_0 \in \widetilde{\mathcal A}_S^{\rm (II)}$.
Thus $A_{f(1)} \cap A_{f(2)}$ is maximal in 
$\{A_{f(1)} \cap A : A\in {\mathcal A}-\{A_{f(1)}\} \}$.

Let $f(3)\in \{ 1,\ldots,l \}-\{f(1),f(2)\}$ 
such that 
$(A_{f(1)} \cup A_{f(2)})\cap A_{f(3)}$ is maximal in 
\[ \{(A_{f(1)} \cup A_{f(2)}) \cap A : A\in \{ A_1,\ldots,A_l \}-\{A_{f(1)},A_{f(2)}\} \}. \]
Then we show that $(A_{f(1)} \cup A_{f(2)})\cap A_{f(3)}$ is maximal in 
\[ \{ (A_{f(1)} \cup A_{f(2)})\cap A : A\in {\mathcal A}-\{A_{f(1)},A_{f(2)}\} \}. \]
Indeed if this does not hold, then 
there exists $A'_3 \in {\mathcal A}-\{A_{f(1)},A_{f(2)} \}$ 
such that $(A_{f(1)} \cup A_{f(2)}) \cap A'_3$ is maximal in 
\[ \{(A_{f(1)} \cup A_{f(2)}) \cap A : A\in {\mathcal A}-\{A_{f(1)},A_{f(2)}\} \} \]
and 
\[ (A_{f(1)} \cup A_{f(2)}) \cap A_{f(3)} \subsetneq (A_{f(1)} \cup A_{f(2)}) \cap A'_3. \]
Here $A'_3 \not\in \{A_1,\ldots,A_l \}$ and $A'_3 \not\subset B_0$.
Hence some $\overline{U}_2 \in \overline{\mathcal U}$ separates $B_0$ and $A'_3$, and
$\overline{U}_2$ separates $A_{f(1)} \cup A_{f(2)}$ and $A'_3$.
Then $\overline{U}_2$ separates $A_{f(1)} \cup A_{f(2)}$ and $A_{f(3)}$, 
because $(A_{f(1)} \cup A_{f(2)}) \cap A_{f(3)} \subsetneq (A_{f(1)} \cup A_{f(2)}) \cap A'_3 \subset \overline{U}_2$.
This contradicts that 
$A_{f(1)} \cup A_{f(2)}\cup A_{f(3)} \subset B_0 \in \widetilde{\mathcal A}_S^{\rm (II)}$.

By iterating this argument, 
we obtain a sequence $A_{f(1)},\ldots,A_{f(l)}$ 
that satisfies the condition $(*)$ in ${\mathcal A}$.
\end{proof}

\begin{Lemma}\label{lemmaA4}
Let $B_1,\ldots,B_n \in \widetilde{\mathcal A}_S$ such that 
the sequence $B_1,\ldots,B_n$  
satisfies the condition $(*)$ in $\widetilde{\mathcal A}_S$.
Then there exists a sequence 
\[ A_1,\ldots,A_{l_1}, A_{l_1+1}, \ldots, A_{l_2}, \ldots, A_{l_{n-1}+1}, \ldots, A_{l_n} \; \in {\mathcal A} \]
satisfying the condition~$(*)$ in ${\mathcal A}$ 
such that 
\[ \textstyle B_1=\bigcup_{i=1}^{l_1} A_i,\ B_2=\bigcup_{i=l_1+1}^{l_2} A_i,\ \ldots,\ 
B_n=\bigcup_{i=l_{n-1}+1}^{l_n} A_i. \]
Also $(B_1 \cup \cdots \cup B_i) \cap B_{i+1}=(B_1 \cup \cdots \cup B_i) \cap A_{l_i+1}$ 
for any $i=1,\ldots,n-1$.
\end{Lemma}

\begin{proof}
In the case that $n=1$, by Lemma~\ref{lemmaA2}, 
there exists a sequence $A_1,\ldots,A_{l_1}\in {\mathcal A}$ satisfying the condition~$(*)$ in ${\mathcal A}$ 
such that $B_1=\bigcup_{i=1}^{l_1} A_i$.

Let $n\ge 2$.
Suppose that $B_1,\ldots,B_{n-1},B_n \in \widetilde{\mathcal A}_S$ 
satisfies the condition $(*)$ in $\widetilde{\mathcal A}_S$, 
a sequence 
\[ A_1,\ldots,A_{l_1}, A_{l_1+1}, \ldots, A_{l_2}, \ldots, A_{l_{n-2}+1}, \ldots, A_{l_{n-1}} \; \in {\mathcal A} \]
satisfies the condition~$(*)$ in ${\mathcal A}$ and 
\[ \textstyle B_1=\bigcup_{i=1}^{l_1} A_i,\ \; B_2=\bigcup_{i=l_1+1}^{l_2} A_i,\ \; \ldots,\ \; 
B_{n-1}=\bigcup_{i=l_{n-2}+1}^{l_{n-1}} A_i. \]
Here $B_n\in \widetilde{\mathcal A}_S$ and 
$(B_1 \cup \cdots \cup B_{n-1}) \cap B_n$ is maximal in 
\[ \{(B_1 \cup \cdots \cup B_{n-1})\cap B : B\in \widetilde{\mathcal A}_S -\{ B_1,\ldots,B_{n-1} \} \}. \]
Since $B_n \not\in \{ B_1,\ldots,B_{n-1} \}$, 
by the definition of $\widetilde{\mathcal A}_S$, 
there exists $\overline{U} \in \overline{\mathcal{U}}$ such that 
$B_n \subset \overline{X}_1$ 
and $B_1 \cup \cdots \cup B_{n-1} \subset \overline{X}_2 \cup \cdots \cup \overline{X}_t$ 
where 
$S-\overline{U}=X_1 \cup \cdots \cup X_t$ is a disjoint union, 
$X_1,\ldots, X_t$ are the connected components of $S-\overline{U}$ and 
$\overline{X}_j=\bigcup \{ A\in {\mathcal A} : A \subset X_j\cup \overline{U} \}$ 
for $j=1,\ldots,t$.
Here $\overline{U}$ is a separator of ${\mathcal A}$.
For $U_0:=\overline{U} \cap \overline{X}_1$, 
by (v) in the definition of a separator, 
there exists $A'_0 \in {\mathcal A}$ such that 
$U_0 \subset A'_0 \subset \overline{X}_1$.
Then by (iv) in the definition of a separator, 
the following equation holds;
\begin{align*}
U_0 &= (\overline{X}_2 \cup \cdots \cup \overline{X}_t)\cap A'_0 =\overline{U} \cap A'_0 \\
&= (\overline{X}_2 \cup \cdots \cup \overline{X}_t)\cap \overline{X}_1 =\overline{U} \cap \overline{X}_1.
\end{align*}

Let $U'_0:=(B_1 \cup \cdots \cup B_{n-1}) \cap B_n$.
Then 
\[ U'_0=(B_1 \cup \cdots \cup B_{n-1}) \cap B_n \subset 
(\overline{X}_2 \cup \cdots \cup \overline{X}_t)\cap \overline{X}_1=U_0. \]
Hence $U'_0\subset U_0 \subset A'_0 \subset \overline{X}_1$.
Thus
\begin{align*}
U'_0 &= (B_1 \cup \cdots \cup B_{n-1}) \cap A'_0 \\
&= (B_1 \cup \cdots \cup B_{n-1}) \cap \overline{X}_1.
\end{align*}
Here $U'_0=(B_1 \cup \cdots \cup B_{n-1}) \cap B_n$ is maximal in 
\[ \{(B_1 \cup \cdots \cup B_{n-1})\cap B : B\in \widetilde{\mathcal A}_S -\{ B_1,\ldots,B_{n-1} \} \}. \]
By the above argument, $U'_0=(B_1 \cup \cdots \cup B_{n-1}) \cap A'_0$ is maximal in 
\[ \{(B_1 \cup \cdots \cup B_{n-1})\cap A : A\in {\mathcal A} -\{ A_1,\ldots,A_{l_{n-1}} \} \} \]
and it is maximal in 
\[ \{(A_1 \cup \cdots \cup A_{l_{n-1}})\cap A : A\in {\mathcal A} -\{ A_1,\ldots,A_{l_{n-1}} \} \}. \]

Hence the sequence 
\[ A_1,\ldots,A_{l_1}, A_{l_1+1}, \ldots, A_{l_2}, \ldots, A_{l_{n-2}+1}, \ldots, A_{l_{n-1}}, \, A'_0 \; \in {\mathcal A} \]
satisfies the condition~$(*)$ in ${\mathcal A}$.
Here $U'_0$ is a separator of ${\mathcal A}$.

Let $A'_1,\ldots,A'_{l'}$ be the set of ${\mathcal A}$ such that $A'_i \subset B_n$ for any $i=1,\ldots,l'$.
By Lemma~\ref{lemmaA1}, $B_n=A'_1 \cup \cdots \cup A'_{l'}$.
By Lemma~\ref{lemmaA2}, we may suppose that $A'_1=A'_0$ and 
the sequence $A'_1,\ldots,A'_{l'}$ satisfies the condition $(*)$ in ${\mathcal A}$.
Then 
\[ (B_1 \cup \cdots \cup B_{n-1}) \cap B_n = U'_0 = (B_1 \cup \cdots \cup B_{n-1}) \cap A'_1. \]

We show that 
\[ ((B_1 \cup \cdots \cup B_{n-1})\cup A'_1 \cup \cdots \cup A'_i)\cap A'_{i+1}=(A'_1 \cup \cdots \cup A'_i)\cap A'_{i+1} \]
for any $i=1,\ldots,l'-1$.
Let $U_{i+1}:=((B_1 \cup \cdots \cup B_{n-1})\cup A'_1 \cup \cdots \cup A'_i)\cap A'_{i+1}$.
Then 
\begin{align*}
U_{i+1}&=((B_1 \cup \cdots \cup B_{n-1})\cup A'_1 \cup \cdots \cup A'_i)\cap A'_{i+1} \\
&\subset (\overline{X}_2 \cup \cdots \cup \overline{X}_t)\cup A'_1 \cup \cdots \cup A'_i)\cap A'_{i+1} \\
&= ((\overline{X}_2 \cup \cdots \cup \overline{X}_t)\cap A'_{i+1}) \cup ((A'_1 \cup \cdots \cup A'_i)\cap A'_{i+1}) \\
&= ((\overline{X}_2 \cup \cdots \cup \overline{X}_t)\cap \overline{X}_1 \cap A'_{i+1}) \cup ((A'_1 \cup \cdots \cup A'_i)\cap A'_{i+1}) \\
&= ((U_0 \cap A'_{i+1}) \cup ((A'_1 \cup \cdots \cup A'_i)\cap A'_{i+1}) \\
&= (U_0 \cup A'_1 \cup \cdots \cup A'_i)\cap A'_{i+1} \\
&= (A'_1 \cup \cdots \cup A'_i)\cap A'_{i+1},
\end{align*}
where $A'_{i+1} \subset \overline{X}_1$ and $A'_{i+1}=\overline{X}_1 \cap A'_{i+1}$.
Also here $(\overline{X}_2 \cup \cdots \cup \overline{X}_t)\cap \overline{X}_1=U_0$ and $U_0 \subset A'_0=A'_1$.
Obviously $(A'_1 \cup \cdots \cup A'_i)\cap A'_{i+1} \subset U_{i+1}$ holds.

Then $U_{i+1}$ is maximal in 
\[ \{(A'_{1}\cup \cdots \cup A'_{i}) \cap A : A \in {\mathcal A}-\{A'_{1},\ldots,A'_{i}\} \} \]
for any $i=1,\ldots,l'-1$.
Hence $U_{i+1}$ is maximal in 
\begin{multline*}
\{((B_1 \cup \cdots \cup B_{n-1})\cup A'_{1}\cup \cdots \cup A'_{i}) \cap A : \\
A \in {\mathcal A}
-\{ A_1,\ldots,A_{l_1}, A_{l_1+1}, \ldots, A_{l_2}, \ldots, A_{l_{n-2}+1}, \ldots, A_{l_{n-1}} ,A'_{1},\ldots,A'_{i}\} \}
\end{multline*}
for any $i=1,\ldots,l'-1$.

Thus the sequence 
\[ A_1,\ldots,A_{l_1}, A_{l_1+1}, \ldots, A_{l_2}, \ldots, A_{l_{n-2}+1}, \ldots, A_{l_{n-1}}, 
A'_1  \ldots, A'_{l'} \; \in {\mathcal A} \]
satisfies the condition~$(*)$ in ${\mathcal A}$.
\end{proof}

\medskip

(4) Now we show that (4) in the definition of a separation of $S$ holds for $\widetilde{\mathcal A}_S$.

Suppose that $B_1,\ldots,B_n \in \widetilde{\mathcal A}_S$ is a sequence 
satisfying the condition~$(*)$ in $\widetilde{\mathcal A}_S$.
Then we show that $U_i:=(B_1\cup \cdots \cup B_i) \cap B_{i+1}$ is a separator of $\widetilde{\mathcal A}_S$ 
for any $i=1,\ldots,n-1$.

Let $i\in \{1,\ldots,n-1\}$ and let $U:=U_i=(B_1\cup \cdots \cup B_i) \cap B_{i+1}$.
By Lemma~\ref{lemmaA4}, $U$ is a separator of ${\mathcal A}$.

\medskip

(i) and (ii): Then $U$ is a spherical-product subset of $S$ and $U$ separates $S$.

\medskip

(iii) 
Suppose that $X-U=X_1 \cup \cdots \cup X_t$ is a disjoint union 
and $X_1,\ldots, X_t$ are the connected components of $S-U$.
Let 
\begin{align*}
&\overline{X}_j:=\bigcup \{ A\in {\mathcal A} : A \subset X_j\cup U \}\ \text{and} \\
&\overline{X}'_j:=\bigcup \{ B\in \widetilde{\mathcal A}_S : B \subset X_j\cup U \} 
\end{align*}
for each $j=1,\ldots,t$.

By the definition of $\widetilde{\mathcal A}_S$, 
each $B\in \widetilde{\mathcal A}_S$ is not separated by any separator of a minimal separation of $S$.
Hence for any $B\in \widetilde{\mathcal A}_S$, $U$ does not separate $B$ and 
$B-U \subset X_{j_0}$ for some unique $j_0\in\{1,\ldots,t\}$.
Then $B \subset \overline{X}_{j_0}$ and $B \subset \overline{X}'_{j_0}$ hold.

Here $\overline{X}_{j}=\overline{X}'_{j}$ holds for any $j=1,\ldots,t$, 
because each $B\in \widetilde{\mathcal A}_S$ is denoted by $B=\bigcup_{i=1}^l A_i$ 
(where $A_1,\ldots,A_l$ are the elements of ${\mathcal A}$ such that $A_i \subset B$) 
and each $A\in {\mathcal A}$ is contained in some unique $B \in \widetilde{\mathcal A}_S$.

\medskip

(iv) 
Since $U$ is a separator of ${\mathcal A}$, 
there exist $A_1,A_2\in {\mathcal A}$ such that 
$A_1 \cap A_2=U$, $A_1 \subset \overline{X}_{j_1}$ and $A_2 \subset \overline{X}_{j_2}$ 
for some $j_1, j_2\in \{1,\ldots,t\}$ as $j_1\neq j_2$.
Then $A_1 \subset B'_1$ and $A_2 \subset B'_2$ 
for some unique $B'_1,B'_2 \in \widetilde{\mathcal A}_S$.
Here $B'_1 \subset \overline{X}_{j_1}=\overline{X}'_{j_1}$ and 
$B'_2 \subset \overline{X}_{j_2}=\overline{X}'_{j_2}$.
Then $B'_1\cap B'_2=U$ holds.

\medskip

(v) 
Let $j\in \{1,\ldots,t\}$.
Suppose that $B'_1,\ldots,B'_k \in \widetilde{\mathcal A}_S$ 
satisfies the condition~$(*)$ in $\widetilde{\mathcal A}_S$ 
and $U \cap \overline{X}'_{j} \subset B'_1\cup \cdots \cup B'_k \subset \overline{X}'_{j}$.
By Lemma~\ref{lemmaA4}, there exists a sequence 
\[ A'_1,\ldots,A'_{p_1}, \ldots, A'_{p_{k-1}+1}, \ldots, A'_{p_k} \; \in {\mathcal A} \]
satisfying the condition~$(*)$ in ${\mathcal A}$ 
such that 
\[ \textstyle B'_1=\bigcup_{i=1}^{p_1} A'_i,\ \ldots,\ B'_k=\bigcup_{i=p_{k-1}+1}^{p_k} A'_i. \]
Since $U$ is a separator of ${\mathcal A}$, 
we have that $U \cap \overline{X}'_{j} \subset A'_{i_0}$ for some $i_0 \in\{1,\ldots,p_k \}$.
Here $A'_{i_0} \subset B'_{i_1}$ for some $i_1\in\{1,\ldots,k \}$.
Then $U \cap \overline{X}'_{j} \subset A'_{i_0} \subset B'_{i_1}$.

\medskip

Thus from the above (i)--(v), $U$ is a separator of $\widetilde{\mathcal A}_S$.

\medskip

We show (b) on (4) in the definition of a separation for $\widetilde{\mathcal A}_S$; that is, 
$B_1\cup \cdots \cup B_i \subset X_{j_1} \cup U$ and $B_{i+1} \subset X_{j_2} \cup U$ 
for some $j_1, j_2\in \{1,\ldots,t \}$ as $j_1 \neq j_2$.

Since $\widetilde{\mathcal A}$ is a separation, 
by Lemma~\ref{lemmaA4}, we obtain that 
$B_1\cup \cdots \cup B_i \subset X_{j_1} \cup U$ and $A_{l_i+1} \subset X_{j_2} \cup U$ 
for some $j_1, j_2\in \{1,\ldots,t \}$ as $j_1 \neq j_2$.
Since $U$ is a separator of $\widetilde{\mathcal A}_S$, 
$U$ does not separate $B_{i+1}=A_{l_i+1} \cup \cdots \cup A_{l_{i+1}} \in \widetilde{\mathcal A}_S$.
Hence $B_{i+1} \subset X_{j_2} \cup U$.

\smallskip

Thus (4) in the definition of a separation of $S$ holds for $\widetilde{\mathcal A}_S$.

\smallskip

Therefore $\widetilde{\mathcal A}_S$ is a separation of $S$.

\section{Type(I)-type(II)-compatible and conjugate up to finite twists}\label{sec6}

We investigate on type(I)-type(II)-compatible Coxeter systems and conjugate up to finite twists.
We show that for Coxeter systems $(W,R)$ and $(W,S)$ with the untangle-condition, 
$R$ and $S$ are conjugate up to finite twists if and only if 
$(W,R)$ and $(W,S)$ are type(I)-type(II)-compatible.

We show a lemma on conjugate spherical-product subsets.

\begin{Lemma}\label{lem0-00}
Let $(W,S)$ be a Coxeter system and 
let $U$ be a spherical-product subset of $S$.
If $wUw^{-1}=U$ for some $w\in W$, then 
$wU_\sigma w^{-1}=U_\sigma$ and $wU_\nu w^{-1}=U_\nu$.
\end{Lemma}

\begin{proof}
Suppose that $wUw^{-1}=U$ for some $w\in W$.
Then 
\begin{align*}
&W_U=W_{U_\sigma} \times W_{U_\nu} \ \text{and}  \\
&W_{wUw^{-1}}=wW_U w^{-1}=wW_{U_\sigma}w^{-1} \times wW_{U_\nu}w^{-1} \\
&\hspace*{14.4mm} =W_{wU_\sigma w^{-1}} \times W_{wU_\nu w^{-1}}.
\end{align*}
Here $W_{U_\nu}$ and $W_{wU_\nu w^{-1}}$ are 
the minimal standard (parabolic) subgroups of finite index 
in $W_U$ and $W_{wUw^{-1}}$ respectively (see \cite{Ho1}).
Since $W_U=W_{wUw^{-1}}$ and $U=wUw^{-1}$, we have that 
$W_{U_\nu}=W_{wU_\nu w^{-1}}$.
Hence $U_\nu=wU_\nu w^{-1}$.
Thus $wU_\nu w^{-1}=U_\nu$ and $wU_\sigma w^{-1}=U_\sigma$.
\end{proof}

We show a lemma on spherical subsets and type(II) subsets.

\begin{Lemma}\label{lem0-2}
Let $(W,S)$ be a Coxeter system.
If $\sigma$ is a spherical subset of $S$ that separates $S$, 
then $\sigma$ does not separate any type(II) subset of $S$.
\end{Lemma}

\begin{proof}
Let $\sigma$ be a spherical subset of $S$ that separates $S$.
We can denote $S - \sigma =X_1 \cup \cdots \cup X_t$ 
where $X_1,\ldots,X_t$ are the connected components of $S - \sigma$.

Then for any maximal twist-rigid subset $A$ of $S$, 
$\sigma$ does not separate $A$ and 
$A \subset X_i \cup \sigma$ for some $i\in \{1,\ldots,t\}$.
Also for any minimal separation ${\mathcal A}$ of $S$ and for any $B \in {\mathcal A}$, 
$\sigma$ does not separate $B$ and 
$B \subset X_i \cup \sigma$ for some $i\in \{1,\ldots,t\}$.
(Indeed if $\sigma$ separates some $B \in {\mathcal A}$, then 
there exists a separation ${\mathcal A}'$ of $S$ induced by $\sigma$ from ${\mathcal A}$ 
such that ${\mathcal A}' \prec {\mathcal A}$.
This contradicts that ${\mathcal A}$ is minimal.)

Hence $\sigma$ does not separate any type(II) subset of $S$.
\end{proof}

Now we show the following.

\begin{Lemma}\label{lem2}
Let $(W,S)$ be a Coxeter system with the untangle-condition 
and let $(W,S')$ be a Coxeter system obtained from $(W,S)$ by some twist.
Then $(W,S)$ and $(W,S')$ are type(I)-type(II)-compatible.
\end{Lemma}

\begin{proof}
It is sufficient to show this lemma 
in the case that $(W,S)$ and $(W,S')$ are connected.
Now we suppose this.

Let $U$ be a spherical-product subset of $S$ and let $w\in W$ 
such that $U$ separates $S$, $U=wUw^{-1}$ and 
$S'$ is obtained from $S$ by some twist of $U$ and $w$.
There exist non-empty subsets $X$ and $Y$ of $S$ 
such that $S-U=X \cup Y$, $X \cap Y=\emptyset$, 
$o(xy)=\infty$ for any $x\in X$ and $y\in Y$, 
and $S'=X \cup U \cup (w Y w^{-1})$.

Let ${\mathcal A}_0$ and ${\mathcal A}'_0$ be the sets of 
maximal twist-rigid subsets of $S$ and $S'$ respectively.
Here for each $C \in {\mathcal A}_0$, 
if we put $C':=C$ (if $C \subset X\cup U$) and 
$C':=w^{-1}Cw$ (if $C\not\subset X\cup U$), 
then $C' \in {\mathcal A}'_0$.
Hence 
\begin{equation*}
\begin{split}
{\mathcal A}'_0
= &\{ C : C\in {\mathcal A}_0,\ C\subset X\cup U \} \\
&\ \hspace*{12mm} 
\cup \{ wCw^{-1} : C\in {\mathcal A}_0,\ C\not\subset X\cup U \}.
\end{split}
\end{equation*}

By the definition of type(I) subsets and Remark~\ref{Rem3-2}, 
each $A\in \widetilde{\mathcal A}_S^{\rm (I)}$ is not separated 
by any spherical-product subset that separates $S$.

Let $A\in \widetilde{\mathcal A}_S^{\rm (I)}$.
Here $A$ is not separated by $U$.
We consider the subset $A'$ of $S'$ as $A'=A$ (if $A\subset X\cup U$) and 
$A'=wAw^{-1}$ (if $A\not\subset X\cup U$).
Then since $U$ separates $S'$, 
there is an induced separation ${\mathcal A}'$ of $S'$ by $U$ as in Remark~\ref{Rem3-2}.
Then $A'$ is not separated by any spherical-product subset $U'$ that separates $S'$.
Indeed, suppose that $A'$ is separated by some spherical-product subset $U'_0$ that separates $S'$.
Here either $U$ does not separate any $B\in \widetilde{\mathcal A}_S^{\rm (II)}$ or 
$U$ separates some $B \in \widetilde{\mathcal A}_S^{\rm (II)}$.
In both cases, $U'_0 \subset X\cup U$ or $U'_0 \subset wYw^{-1}\cup U$.
For $U_0:=U'_0$ (if $U'_0 \subset X\cup U$) and 
$U_0:=w^{-1}U'_0w$ (if $U'_0\not\subset X\cup U$), 
$U_0$ is a spherical-product subset of $S$ that separates $S$ and $U_0$ has to separate $A$ 
that is a contradiction.
Hence $A'\in \widetilde{\mathcal A}_{S'}^{\rm (I)}$.

Also by the same argument, 
for each $A'\in \widetilde{\mathcal A}_{S'}^{\rm (I)}$, 
if $A$ is the subset of $S$ defined by 
$A:=A'$ (if $A'\subset X\cup U$) and 
$A:=w^{-1}A'w$ (if $A'\not\subset X\cup U$) 
then $A\in \widetilde{\mathcal A}_{S}^{\rm (I)}$.

Thus 
\begin{equation*}
\begin{split}
\widetilde{\mathcal A}_{S'}^{\rm (I)} 
= &\{ A : A\in \widetilde{\mathcal A}_S^{\rm (I)},\ A\subset X\cup U \} \\
&\ \hspace*{12mm} 
\cup \{ wAw^{-1} : A\in \widetilde{\mathcal A}_S^{\rm (I)},\ A\not\subset X\cup U \},
\end{split}
\end{equation*}
and $|\widetilde{\mathcal A}_S^{\rm (I)}|=|\widetilde{\mathcal A}_{S'}^{\rm (I)}|$.

\medskip

(a) We first suppose that 
$U$ does not separate any $B\in \widetilde{\mathcal A}_S^{\rm (II)}$.

Then for any minimal separation ${\mathcal A}$ of $S$, 
$S'$ is a twist of $S$ that is preserving ${\mathcal A}$ and 
the separation ${\mathcal A}'$ of $S'$ induced by 
${\mathcal A}$ and the twist is minimal.
Also for any minimal separation ${\mathcal A}'$ of $S'$, 
the separation ${\mathcal A}$ of $S$ induced by 
${\mathcal A}'$ and the twist is minimal.
Hence 
\begin{equation*}
\begin{split}
\widetilde{\mathcal A}_{S'}^{\rm (II)} 
= &\{ B : B\in \widetilde{\mathcal A}_S^{\rm (II)},\ B\subset X\cup U \} \\
&\ \hspace*{12mm} 
\cup \{ wBw^{-1} : B\in \widetilde{\mathcal A}_S^{\rm (II)},\ B\not\subset X\cup U \}.
\end{split}
\end{equation*}
Thus $(W,S)$ and $(W,S')$ are type(I)-type(II)-compatible.

\medskip

(b) We suppose that 
$U$ separates some unique $B_0\in \widetilde{\mathcal A}_S^{\rm (II)}$.
Here $U$ does not separate any $A\in \widetilde{\mathcal A}_S-\{B_0\}$.

Let 
\begin{align*}
&\overline{\mathcal U}_S:=
\bigcap\{{\mathcal U}_{\mathcal A}:
\text{${\mathcal A}$ is a minimal separation of $S$} \} \ \text{and} \\
&\overline{\mathcal U}_{S'}:=
\bigcap\{{\mathcal U}_{{\mathcal A}'}:
\text{${\mathcal A}'$ is a minimal separation of $S'$} \},
\end{align*}
where ${\mathcal U}_{\mathcal A}$ and ${\mathcal U}_{{\mathcal A}'}$ 
are the sets of separators of 
minimal separations ${\mathcal A}$ and ${\mathcal A}'$ of $S$ and $S'$ respectively.
The standard separations $\widetilde{\mathcal A}_S$ and $\widetilde{\mathcal A}_{S'}$ are defined 
by the separator sets $\overline{\mathcal U}_S$ and $\overline{\mathcal U}_{S'}$ 
respectively.

Then for any minimal separation ${\mathcal A}$ of $S$, 
$S'$ is a twist of $S$ that is preserving $\{ A\in {\mathcal A}:A \not\subset B_0 \}$; 
that is, $A \subset X\cup U$ or $A \subset wYw^{-1}\cup U$ for each $A\in {\mathcal A}$ as $A \not\subset B_0$.
Also for any minimal separation ${\mathcal A}'$ of $S'$, 
$S$ is a twist of $S'$ that is preserving $\{ A'\in {\mathcal A}':A' \not\subset B'_0 \}$ 
where $B'_0:=(B_0\cap (X\cup U))\cup w(B_0\cap Y)w^{-1}$; 
that is, $A' \subset X\cup U$ or $A' \subset w^{-1}Yw \cup U$ for each $A'\in {\mathcal A}'$ as $A' \not\subset B'_0$.

Hence for any $\overline{U} \in {\mathcal U}_{\mathcal A}$ as $\overline{U}\subset X \cup U$, 
$\overline{U} \in {\mathcal U}_{{\mathcal A}'}$.
Also for any $\overline{U} \in {\mathcal U}_{\mathcal A}$ as $\overline{U}\subset Y \cup U$, 
$w \overline{U}w^{-1} \in {\mathcal U}_{{\mathcal A}'}$.

Let $B\in \widetilde{\mathcal A}_S^{\rm (II)}-\{ B_0 \}$.
We consider the subset $B'$ of $S'$ as 
$B'=B$ (if $B\subset X\cup U$) and 
$B'=wBw^{-1}$ (if $B\not\subset X\cup U$).
Then $B'\in \widetilde{\mathcal A}_{S'}^{\rm (II)}$ by the above argument.
Here $B$ and $B'$ are type(II)-compatible, 
since $B$ and $B'$ are conjugate.
Hence 
\begin{equation*}
\begin{split}
\widetilde{\mathcal A}_{S'}^{\rm (II)} 
\supset &\{ B : B\in \widetilde{\mathcal A}_S^{\rm (II)}-\{B_0\},\ B\subset X\cup U \} \\
&\ \hspace*{12mm} 
\cup \{ wBw^{-1} : B\in \widetilde{\mathcal A}_S^{\rm (II)}-\{B_0\},\ B\not\subset X\cup U \}.
\end{split}
\end{equation*}
Thus $|\widetilde{\mathcal A}_{S'}^{\rm (II)}|\ge |\widetilde{\mathcal A}_S^{\rm (II)}|$, 
because $|\widetilde{\mathcal A}_S^{\rm (II)}-\{B_0\}|=|\widetilde{\mathcal A}_S^{\rm (II)}|-1$ 
and there is a possibility that 
the corresponding subset $B'_0=(B_0\cap (X\cup U))\cup w(B_0\cap Y)w^{-1}$ of $S'$ to $B_0$ in $S$ is 
a union of some type(II) subsets.
Also since $S$ is a twist of $S'$, by the same argument, 
$|\widetilde{\mathcal A}_S^{\rm (II)}|\ge |\widetilde{\mathcal A}_{S'}^{\rm (II)}|$.
Hence 
$|\widetilde{\mathcal A}_S^{\rm (II)}| = |\widetilde{\mathcal A}_{S'}^{\rm (II)}|$.
This implies that 
\[ B'_0=(B_0\cap (X\cup U))\cup w(B_0\cap Y)w^{-1} \]
is a type(II) subset of $S'$.
Here $B'_0$ is a twist of $B_0$ that induces the twist of $S$.

Thus, if $U$ separates some unique $B_0\in \widetilde{\mathcal A}_S^{\rm (II)}$ 
then $(W,S)$ and $(W,S')$ are type(I)-type(II)-compatible.

\medskip

(c) We suppose that 
$U$ separates some $B\in \widetilde{\mathcal A}_S^{\rm (II)}$ that is not unique.
Let $B_1,\ldots, B_l$ be the type(II) subsets of $S$ 
such that $U$ (also $U\cap B_i$) separates $B_i$ for $i=1,\ldots,l$.

For each $i=1,\ldots,l$, 
we show that $U\cap B_i$ separates unique $B_i\in \widetilde{\mathcal A}_S^{\rm (II)}$.
We suppose that $U\cap B_i$ separates some $B_j$ as $j\neq i$.
Then $U\cap B_i \cap B_j$ separates $B_j$.
The standard separation $\widetilde{\mathcal A}_S$ of $S$ is defined 
by the separator-set $\overline{\mathcal U}_S$.
Here $B_i \cap B_j \subset U_0$ for some $U_0 \in \overline{\mathcal U}_S$.
Then $U\cap B_i \cap B_j \subset U \cap U_0 \subset U_0$.
Since $U\cap B_i \cap B_j$ separates $B_j$, we have that $U_0$ separates $B_j$.
This contradicts that $U_0 \in \overline{\mathcal U}_S$ does not separate 
any $B \in \widetilde{\mathcal A}_S$.

Now $S$ is connected.
Let $Y_1,\ldots, Y_t$ be the connected components of $Y$.
Then $Y=Y_1\cup\cdots \cup Y_t$ and 
$S-U=X\cup Y=X\cup Y_1\cup\cdots \cup Y_t$ that are disjoint unions.

Let ${\mathcal A}_0$ be the maximal twist-rigid subset of $S$, 
let ${\mathcal A}$ be the induced separation of $S$ by $U$ as Remark~\ref{Rem3-2}, 
and let ${\mathcal A}'$ be a minimal separation of $S$ 
such that ${\mathcal A}' \preceq {\mathcal A}$.

Let $i\in \{1,\ldots,t\}$ and let 
\[ \overline{Y}_i:=\bigcup \{ A\in {\mathcal A}_0 : A \subset Y_i\cup U 
\ \text{and} \ A\not\subset U \}. \]
Then 
\[ \overline{Y}_i=\bigcup \{ A\in {\mathcal A}' : A \subset Y_i\cup U 
\ \text{and} \ A\not\subset U \}. \]
Since ${\mathcal A}'$ is a separation of $S$, 
there exists $A'_i \in {\mathcal A}'$ such that 
$\overline{Y}_i \cap A'_i$ is a separator of ${\mathcal A}'$ 
and $A'_i \subset B_{j_i}$ for some $j_i\in \{1,\ldots,l\}$.
Then $\overline{Y}_i \cap A'_i \subset U \cap B_{j_i}$ separates $B_{j_i}$.
(Also $\overline{Y}_i \cap A'_i \subset U \cap B_{j_i}$ separates $S$.)
Here $B_{j_i}\in \widetilde{\mathcal A}_S^{\rm (II)}$ is uniquely determined.

Hence $t=l$ and the map $g(i)=j_i$ is bijective on the set $\{1,\ldots,l\}$.
We may suppose that $g(i)=i$ for any $i=1,\ldots,l$.
Then for each $i=1,\ldots,l$, 
$\overline{Y}_i \cap A'_i$ is a separator of ${\mathcal A}'$ and $A'_i \subset B_i$.

\smallskip

Let $U_i:=(U\cap B_i)\cup U_\sigma$ for each $i=1,\ldots,l$.
Since $U \cap B_i$ separates $S$, $U_i$ separates $S$.
Here $U_i=U_\sigma \cup (U_\nu\cap U_i)$ holds.
Then $wU_\sigma w^{-1}=U_\sigma$ and $wU_\nu w^{-1}=U_\nu$ by Lemma~\ref{lem0-00}.
Hence $w a w^{-1}=a$ for any $a\in U_\nu$ by the untangle-condition.
Thus 
\begin{align*}
wU_iw^{-1} &=(w U_\sigma w^{-1}) \cup (w(U_\nu\cap U_i) w^{-1}) \\
&=U_\sigma \cup (U_\nu\cap U_i) =U_i,
\end{align*}
for any $i=1,\ldots,l$.

Here $U_i$ separates $B_i\in \widetilde{\mathcal A}_S^{\rm (II)}$.
Indeed
\[ B_i-U_i = B_i-(U_i\cap B_i)=B_i-(U \cap B_i)=B_i -U. \]

We show that $U_i$ separates 
the unique element $B_i$ of $\widetilde{\mathcal A}_S^{\rm (II)}$.
We suppose that $U_i$ separates $B_j \in \widetilde{\mathcal A}_S^{\rm (II)}$ as $j\neq i$.
Then 
\[ B_j \cap U_i=(U\cap B_i\cap B_j)\cup (U_{\sigma} \cap B_j) \]
separates $B_j$.
Since the spherical subset $U_{\sigma} \cap B_j$ does not separate 
the type(II) subset $B_j$ by Lemma~\ref{lem0-2}, 
$U\cap B_i\cap B_j \neq \emptyset$ and $B_i\cap B_j \neq \emptyset$.
Hence $B_i\cup B_j$ is connected.
There exists the separation ${\mathcal A}_1$ induced by $U_i$ as in Remark~\ref{Rem3-2} and 
we can obtain a minimal separation ${\mathcal A}'_1$ as ${\mathcal A}'_1 \preceq {\mathcal A}_1$.
Here the spherical-product subset $U_i=(U\cap B_i)\cup U_\sigma$ separates $B_i$ and 
$U_i$ separates $B_j$ by hypothesis.
This contradicts that $B_i$ and $B_j$ are distinct type(II) subsets of $S$.
Thus $U_i$ separates unique $B_i \in \widetilde{\mathcal A}_S^{\rm (II)}$.

\smallskip

For each $i=1,\ldots,l$ and for $X_i:=S-(U_i\cup Y_i)$, 
$S-U_i=X_i \cup Y_i$ that is a disjoint union 
and $o(xy)=\infty$ for any $x\in X_i$ and $y\in Y_i$.
Here 
\[ S=X \cup U \cup (Y_1\cup\cdots \cup Y_l) \]
is a disjoint union and $wU_iw^{-1}=U_i$ as above.

Let $S_1:=X_1 \cup U_1 \cup wY_1w^{-1}$ that is a twist of $S$ by $U_1$ and $w$.
Then $U_1$ separates unique $B_1\in \widetilde{\mathcal A}_S^{\rm (II)}$ by the above argument.
Hence $S$ and $S_1$ are type(I)-type(II)-compatible by the above argument (b).
Then 
\[ S_1=X \cup U \cup (wY_1w^{-1}\cup Y_2 \cup Y_3\cup\cdots \cup Y_l). \]
For the corresponding subset $B'_2$ of $S_1$ to the subset $B_2$ of $S$, 
$U_2$ separates unique $B'_2\in \widetilde{\mathcal A}_{S_1}^{\rm (II)}$ in $S_1$.
Also for the corresponding subset $X'_2$ of $S_1$ to the subset $X_2$ of $S$, 
$S_1-U_2=X'_2 \cup Y_2$ is a disjoint union.

Let $S_2:=X'_2 \cup U_2 \cup wY_2w^{-1}$ that is a twist of $S_1$ by $U_2$ and $w$.
Then $U_2$ separates unique $B'_2\in \widetilde{\mathcal A}_{S_2}^{\rm (II)}$ in $S_2$.
Hence $S_1$ and $S_2$ are type(I)-type(II)-compatible by the above argument (b).
Here
\[ S_2=X \cup U \cup (wY_1w^{-1}\cup wY_2w^{-1} \cup Y_3\cup\cdots \cup Y_l). \]

By iterating this argument,
we obtain a sequence of Coxeter generating sets 
$S=S_0,S_1, S_2,\ldots, S_l$ such that 
$S_i$ is some twist of $S_{i-1}$ by $U_i$ and $w$ for each $i=1,\ldots,l$.
Here $U_i$ separates unique $B'_i\in \widetilde{\mathcal A}_{S_i}^{\rm (II)}$ 
in $(W,S_i)$ for any $i=1,\ldots,l$, 
where $B'_i$ is the type(II) subset of $S_i$ corresponding to 
$B_i$ in $S$.
Then $(W,S_i)$ and $(W,S_{i+1})$ are type(I)-type(II)-compatible for each $i=1,\ldots,l$ 
by the above argument (b).
Here
\begin{align*}
&S=S_0=X \cup U \cup (Y_1\cup Y_2\cup Y_3\cup\cdots \cup Y_l), \\
&S_1=X \cup U \cup (wY_1w^{-1}\cup Y_2 \cup Y_3\cup\cdots \cup Y_l), \\
&S_2=X \cup U \cup (wY_1w^{-1}\cup wY_2w^{-1} \cup Y_3\cup\cdots \cup Y_l), \\
&\ \cdots \\
&S_l=X \cup U \cup (wY_1w^{-1}\cup wY_2w^{-1} \cup\cdots \cup wY_lw^{-1})=S'.
\end{align*}
Thus, we obtain that $(W,S)$ and $(W,S')$ are type(I)-type(II)-compatible.
\end{proof}

\begin{Remark}\label{lem2-0}
Let $(W,S)$ be a Coxeter system.
Let $U_1$ and $U_2$ be spherical-product subsets of $S$ that separate $S$ 
such that for each $i=1,2$, 
\begin{enumerate}
\item[$\text{(a)}_i$] $U_i$ does not separate any $B\in \widetilde{\mathcal A}_{S}^{\rm (II)}$, or
\item[$\text{(b)}_i$] $U_i$ separates some unique $B\in \widetilde{\mathcal A}_{S}^{\rm (II)}$.
\end{enumerate}
We consider a twist $S'=X\cup U_1 \cup wYw^{-1}$ of $S$ by $U_1$ and $w\in W$.
Here we suppose that $wU_1w^{-1}=U_1$, $S-U_1=X \cup Y$ is a disjoint union and 
$o(xy)=\infty$ for any $x\in X$ and $y\in Y$.

Then we investigate a spherical-product subset $U'_2$ of $S'$ corresponding to $U_2$ in $S$ 
from the proof of Lemma~\ref{lem2}.

If 
\begin{enumerate}
\item[$\text{(a)}_1$] $U_1$ does not separate any $B\in \widetilde{\mathcal A}_{S}^{\rm (II)}$,
\end{enumerate}
then for $U'_2:=U_2$ (if $U_2 \subset X \cup U_1$) and $U'_2:=w U_2 w^{-1}$ (if $U_2 \not\subset X \cup U_1$), 
$U'_2$ is the spherical-product subset of $S'$ corresponding to $U_2$ in $S$.

We suppose that 
\begin{enumerate}
\item[$\text{(b)}_1$] $U_1$ separates some unique $B_1\in \widetilde{\mathcal A}_{S}^{\rm (II)}$.
\end{enumerate}
If 
\begin{enumerate}
\item[$\text{(a)}_2$] $U_2$ does not separate any $B\in \widetilde{\mathcal A}_{S}^{\rm (II)}$,
\end{enumerate}
then for $U'_2:=U_2$ (if $U_2 \subset X \cup U_1$) and $U'_2:=w U_2 w^{-1}$ (if $U_2 \not\subset X \cup U_1$), 
$U'_2$ is the spherical-product subset of $S'$ corresponding to $U_2$ in $S$.
Also if 
\begin{enumerate}
\item[$\text{(b)}_2$] $U_2$ separates some unique $B_2\in \widetilde{\mathcal A}_{S}^{\rm (II)}$
\end{enumerate}
and if $B_1 \neq B_2$, 
then for $U'_2:=U_2$ (if $U_2 \subset X \cup U_1$) and $U'_2:=w U_2 w^{-1}$ (if $U_2 \not\subset X \cup U_1$), 
$U'_2$ is the spherical-product subset of $S'$ corresponding to $U_2$ in $S$.

Suppose that 
\begin{enumerate}
\item[$\text{(b)}_1$] $U_1$ separates some unique $B_1\in \widetilde{\mathcal A}_{S}^{\rm (II)}$ and
\item[$\text{(b)}_2$] $U_2$ separates some unique $B_2\in \widetilde{\mathcal A}_{S}^{\rm (II)}$
\end{enumerate}
and suppose that $B_1 = B_2$.
In this case, $U_1$ and $U_2$ both separate the unique type(II) subset $B_1=B_2$.
\end{Remark}

We obtain the following theorem from the proof of Lemma~\ref{lem2} and Remark~\ref{lem2-0}.

\begin{Theorem}\label{lem2-1}
Let $(W,R)$ be a Coxeter system with the untangle-condition 
and let $S$ be a Coxeter generating set for $W$ obtained from $R$ by some finite twists.
Then $(W,R)$ and $(W,S)$ are type(I)-type(II)-compatible.
\end{Theorem}

\begin{proof}
There exists a sequence 
$R=S_1,S_2,\ldots,S_n=S$ of Coxeter generating sets for $W$ 
such that $S_{i+1}$ is obtained from $S_i$ by some twist for any $i=1,\ldots,n-1$.
We suppose that the twist $S_{i+1}$ of $S_i$ is obtained 
by a spherical-product subset $U_i$ of $S_i$ and $w_i \in W$ as $w_i U_i w_i^{-1}=U_i$ 
for each $i=1,\ldots,n-1$.
Here by the proof of Lemma~\ref{lem2}, for each $i=1,\ldots,n-1$, either 
\begin{enumerate}
\item[(a)] $U_i$ does not separate any $B\in \widetilde{\mathcal A}_{S_i}^{\rm (II)}$,
\item[(b)] $U_i$ separates some unique $B\in \widetilde{\mathcal A}_{S_i}^{\rm (II)}$, or
\item[(c)] $U_i$ separates some $B\in \widetilde{\mathcal A}_{S_i}^{\rm (II)}$ that is not unique.
\end{enumerate}
Here if (c) $U_i$ separates some $B\in \widetilde{\mathcal A}_{S_i}^{\rm (II)}$ that is not unique, 
then by the proof of Lemma~\ref{lem2}, 
there exist a sequence $U_i^{1},\ldots,U_i^{t}$ of spherical-product subsets 
and a sequence $S_i=S_i^1,S_i^2,\ldots,S_i^t=S_{i+1}$ of Coxeter generating set for $W$ 
such that each $S_i^{j+1}$ is obtained from $S_i^{j}$ by some twist induced by $U_i^{j}$ and $w_i$ 
as $w_i U_i^{j} w_i^{-1}=U_i^{j}$ and 
$U_i^{j}$ separates some unique $B_i^{j}\in \widetilde{\mathcal A}_{S_i^j}^{\rm (II)}$.

Thus there exists a sequence 
$R=S_1,S_2,\ldots,S_k=S$ of Coxeter generating sets for $W$ 
such that for any $i=1,\ldots,k-1$, 
$S_{i+1}$ is obtained from $S_i$ by some twist 
of some spherical-product subset $U_i$ of $S_i$ and $w_i \in W$ as $w_i U_i w_i^{-1}=U_i$ 
and either 
\begin{enumerate}
\item[(a)] $U_i$ does not separate any $B\in \widetilde{\mathcal A}_{S_i}^{\rm (II)}$, or
\item[(b)] $U_i$ separates some unique $B\in \widetilde{\mathcal A}_{S_i}^{\rm (II)}$.
\end{enumerate}

By the proof of Lemma~\ref{lem2} and Remark~\ref{lem2-0},
we can obtain that $(W,R)=(W,S_1)$ and $(W,S)=(W,S_k)$ are type(I)-type(II)-compatible.
\end{proof}

We show the following theorem.

\begin{Theorem}\label{propB}
Let $(W,R)$ and $(W,S)$ be Coxeter systems with the untangle-condition.
If $(W,R)$ and $(W,S)$ are type(I)-type(II)-compatible, 
then there exists a Coxeter generating set $\overline{R}$ 
obtained from $R$ by some finite twists 
such that $(W,\overline{R})$ and $(W,S)$ are compatible on 
the standard separations $\widetilde{\mathcal A}_{\overline{R}}$ and $\widetilde{\mathcal A}_S$ 
(hence they are some-separation-compatible).
\end{Theorem}

\begin{proof}
Let $A_1 \in \widetilde{\mathcal A}_{R}^{\rm (II)}$ and 
$B_1 \in \widetilde{\mathcal A}_S^{\rm (II)}$ such that 
$A_1$ and $B_1$ are type(II)-compatible.
There exists $A'_1$ obtained from $A_1$ by some finite twists 
that induce some twists of $R$ preserving $\widetilde{\mathcal A}_R-\{ A_1 \}$ 
such that $A'_1$ and $B_1$ are conjugate in $W$.
Let $R_1$ be the Coxeter generating set for $W$ obtained 
from $R$ by the induced finite twists preserving $\widetilde{\mathcal A}_R-\{ A_1 \}$ 
such that $A'_1$ in $R_1$ is conjugate to $B_1$ in $S$.

Then by the proofs of Lemma~\ref{lem2} and Theorem~\ref{lem2-1}, 
$(W,R)$ and $(W,R_1)$ are type(I)-type(II)-compatible, 
$A'_1$ is a type(II) subset of $R_1$ and 
\[ \widetilde{\mathcal A}_{R_1} = (\widetilde{\mathcal A}_R - \{ A_1 \})' \cup \{ A'_1 \}, \]
where $(\widetilde{\mathcal A}_R - \{ A_1 \})'$ is the set of the corresponding conjugate subsets of $R_1$ 
to the elements of $\widetilde{\mathcal A}_R - \{ A_1 \}$ in $R$.

Hence 
\begin{enumerate}
\item[(1)] $(W,R_1)$ and $(W,S)$ are type(I)-type(II)-compatible, 
\item[(2)] $A'_1 \in \widetilde{\mathcal A}_{R_1}^{\rm (II)}$ is conjugate to 
$B_1 \in \widetilde{\mathcal A}_S^{\rm (II)}$ and
\item[(3)] $\widetilde{\mathcal A}_{R_1}-\{A'_1\}=(\widetilde{\mathcal A}_R - \{ A_1 \})'$.
\end{enumerate}

Let $A'_2 \in \widetilde{\mathcal A}_{R_1}^{\rm (II)}-\{ A_1'\}$ and 
$B_2 \in \widetilde{\mathcal A}_S^{\rm (II)}-\{ B_1\}$ such that 
$A'_2$ and $B_2$ are type(II)-compatible.
Let $A_2 \in \widetilde{\mathcal A}_{R}^{\rm (II)}$ be the corresponding conjugate subset of $R$ 
to $A'_2 \in \widetilde{\mathcal A}_{R_1}^{\rm (II)}$ in $R_1$.

There exists $A''_2$ obtained from $A'_2$ by some finite twists 
that induce some twists of $R_1$ preserving $\widetilde{\mathcal A}_{R_1}-\{ A'_2 \}$ 
such that $A''_2$ and $B_2$ are conjugate.
Let $R_2$ be the induced Coxeter generating set.

Then by the proofs of Lemma~\ref{lem2} and Theorem~\ref{lem2-1}, 
$(W,R_1)$ and $(W,R_2)$ are type(I)-type(II)-compatible, 
$A''_2$ is a type(II) subset of $R_2$ and 
\[ \widetilde{\mathcal A}_{R_2} = (\widetilde{\mathcal A}_{R_1} - \{ A'_2 \})' \cup \{ A''_2 \} 
= (\widetilde{\mathcal A}_{R} - \{ A_1,\, A_2 \})'' \cup \{ A''_1,\, A''_2 \}, \]
where $A''_1$ is the corresponding conjugate subset of $R_2$ to $A'_1$ in $R_1$, and 
$(\widetilde{\mathcal A}_{R_1} - \{ A'_2 \})'$ (and $(\widetilde{\mathcal A}_{R} - \{ A_1,\, A_2 \})''$) 
is the set of the corresponding conjugate subsets of $R_2$ 
to the elements of $\widetilde{\mathcal A}_{R_1} - \{ A'_2 \}$ in $R_1$ 
(and $\widetilde{\mathcal A}_{R} - \{ A_1,\, A_2 \}$ in $R$ respectively).

Hence
\begin{enumerate}
\item[(1)] $(W,R_2)$ and $(W,S)$ are type(I)-type(II)-compatible, 
\item[(2)] $A''_1 \in \widetilde{\mathcal A}_{R_2}^{\rm (II)}$ is conjugate to 
$B_1 \in \widetilde{\mathcal A}_S^{\rm (II)}$,
\item[(3)] $A''_2 \in \widetilde{\mathcal A}_{R_2}^{\rm (II)}$ is conjugate to 
$B_2 \in \widetilde{\mathcal A}_S^{\rm (II)}$ and
\item[(4)] $\widetilde{\mathcal A}_{R_2}-\{A''_1,\, A''_2\}=(\widetilde{\mathcal A}_R - \{ A_1,\, A_2 \})''$.
\end{enumerate}

By iterating this argument,
we obtain a Coxeter generating set $\overline{R}$ for $W$ 
from $R$ by some finite twists preserving $\widetilde{\mathcal A}_{R}^{\rm (I)}$ 
such that 
for each $A \in \widetilde{\mathcal A}_{R}^{\rm (II)}$, 
the corresponding subset $A' \in \widetilde{\mathcal A}_{\overline{R}}^{\rm (II)}$ 
is conjugate to some unique $B \in \widetilde{\mathcal A}_S^{\rm (II)}$.

Also by the definition of type(I)-type(II)-compatible and the above argument, 
each $A\in \widetilde{\mathcal A}_{\overline{R}}^{\rm (I)}$ 
is conjugate to some unique $B\in \widetilde{\mathcal A}_S^{\rm (I)}$.

Thus, $(W,\overline{R})$ and $(W,S)$ are compatible 
on the standard separations $\widetilde{\mathcal A}_{\overline{R}}$ and $\widetilde{\mathcal A}_{S}$.
\end{proof}

We obtain the following theorem 
from Theorems~\ref{maintheorem}, \ref{lem2-1} and \ref{propB}.

\begin{Theorem}\label{maintheorem2}
Let $(W,R)$ and $(W,S)$ be Coxeter systems with the untangle-condition.
Then the following two statements are equivalent$:$
\begin{enumerate}
\item[\textnormal{(i)}] $R$ and $S$ are conjugate up to finite twists.
\item[\textnormal{(ii)}] $(W,R)$ and $(W,S)$ are type(I)-type(II)-compatible.
\end{enumerate}
\end{Theorem}

\begin{proof}
By Theorem~\ref{lem2-1}, if $R$ and $S$ are conjugate up to finite twists, 
then $(W,R)$ and $(W,S)$ are type(I)-type(II)-compatible.

If $(W,R)$ and $(W,S)$ are type(I)-type(II)-compatible, 
then there exists a Coxeter generating set $\overline{R}$ 
obtained from $R$ by some finite twists 
such that $(W,\overline{R})$ and $(W,S)$ are compatible on 
the standard separations $\widetilde{\mathcal A}_{\overline{R}}$ and $\widetilde{\mathcal A}_S$ 
by Theorem~\ref{propB}.
Thus by Theorem~\ref{maintheorem}, $R$ and $S$ are conjugate up to finite twists.
\end{proof}

%

%
\end{document}